\documentclass[aap,MSNbibl,seceqn,citesort,dvips]{arximspdf}
\usepackage{graphicx}

\doi{10.1214/11-AAP828} 
\volume{22}
\issue{6}
\pubyear{2012}
\firstpage{2320}
\lastpage{2356}

\makeatletter
\newcommand{\eqref}[1]{(\ref{#1})}
\newcommand{\dist}{\stackrel{\mathcal{D}}{\sim}}
\newcommand{\iid}{\mathrm{i.i.d.}}
\newcommand{\tr}{\operatorname{Tr}}
\newcommand{\var}{\operatorname{Var}}
\newcommand{\C}{\mathcal{C}}
\newcommand{\OO}{\mathcal{O}}
\newcommand{\RR}{\mathbb{R}}
\newcommand{\EE}{\mathbb{E}}
\newcommand{\FF}{\mathcal{F}}
\newcommand{\NN}{\mathbb{N}}
\newcommand{\eps}{\varepsilon}
\newcommand{\pphi}{\hat{\varphi}}
\newcommand{\Normal}{\mathrm{N}}
\newcommand{\err}{\mathbf{e}}
\newcommand{\h}{\mathcal{H}}
\newcommand{\eqdef}{\stackrel{\mathrm{def}}{=}}
\newcommand{\weak}{\Rightarrow}
\newcommand{\longweak}{\Longrightarrow}
\newcommand{\opnm}[3]{\|#1\|_{{\mathcal{L}}({\mathcal{H}}^{#2},{\mathcal{H}}^{#3})}}
\renewcommand\phi{\varphi}

\newtheorem{theorem}{Theorem}[section]
\newtheorem{lemma}[theorem]{Lemma}
\newtheorem{corollary}[theorem]{Corollary}
\newtheorem{prop}[theorem]{Proposition}

\newproclaim{remark}[theorem]{Remark}
\newproclaim{remarks}[theorem]{Remarks}
\newproclaim{assumptions}[theorem]{Assumptions}
\newproclaim{assumption}[theorem]{Assumption}
\newproclaim{heuristic}[theorem]{Heuristic}

\setattribute{abstract}{width}{285pt}

\makeatother

\begin{document}
\begin{frontmatter}

\title{Optimal scaling and diffusion limits for the Langevin algorithm in high dimensions}
\runtitle{Langevin algorithm in high dimensions}

\begin{aug}
\author[A]{\fnms{Natesh S.} \snm{Pillai}\thanksref{t1}\ead[label=e1]{pillai@fas.harvard.edu}},
\author[B]{\fnms{Andrew M.} \snm{Stuart}\thanksref{t2}\ead[label=e2]{a.m.stuart@warwick.ac.uk}}
and
\author[C]{\fnms{Alexandre~H.}~\snm{Thi{\'e}ry}\corref{}\thanksref{t3}\ead[label=e3]{a.h.thiery@warwick.ac.uk}}
\thankstext{t1}{Supported by NSF Grant DMS-11-07070.}
\thankstext{t2}{Supported by EPSRC and ERC.}
\thankstext{t3}{Supported by (EPSRC-funded) CRISM.}
\runauthor{N. S. Pillai, A. M. Stuart and A. H. Thi{\'e}ry}
\affiliation{Harvard University, Warwick University and Warwick University}
\address[A]{N. S. Pillai\\
Department of Statistics\\
Harvard University\\
Cambridge, Massachusetts 02138-2901\\
USA\\
\printead{e1}} 
\address[B]{A. M. Stuart\\
Mathematics Institute\\
Warwick University\\
CV4 7AL, Coventry\\
United Kingdom\\
\printead{e2}}
\address[C]{A. H. Thi{\'e}ry\\
Department of Statistics\\
Warwick University\\
CV4 7AL, Coventry\\
United Kingdom\\
\printead{e3}}
\end{aug}

\received{\smonth{3} \syear{2011}}
\revised{\smonth{11} \syear{2011}}

%
\begin{abstract}
The Metropolis-adjusted
Langevin (MALA) algorithm is a sampling algorithm which makes
local moves by incorporating information about the gradient
of the logarithm of the target density.
In this paper we study the efficiency of MALA
on a natural class of target measures supported on an infinite
dimensional Hilbert space. These natural measures have density
with respect to a Gaussian random field measure and arise in many
applications such as
Bayesian nonparametric statistics and the theory of conditioned
diffusions. We prove that, started in stationarity, a suitably
interpolated and scaled version of the Markov chain corresponding
to MALA converges to an infinite dimensional diffusion process.
Our results imply that, in stationarity, the MALA algorithm
applied to an $N$-dimensional approximation of
the target will take $\mathcal{O}(N^{1/3})$ steps
to explore the invariant measure, comparing favorably
with the Random Walk Metropolis which was recently shown to require
$\mathcal{O}(N)$ steps when applied to the same class of problems.
As a by-product of the diffusion limit, it also follows that the
MALA algorithm is optimized at an average acceptance
probability of $0.574$. Previous results were proved
only for targets
which are products of one-dimensional distributions,
or for variants of this situation, limiting their applicability.
The correlation in our target means that the rescaled MALA algorithm converges
weakly to an infinite dimensional Hilbert space valued diffusion,
and the limit cannot be described through analysis of scalar diffusions.
The limit theorem is proved by showing that a
drift-martingale decomposition of the
Markov chain, suitably scaled, closely resembles
a weak Euler--Maruyama discretization of the putative limit.
An invariance principle is proved for the martingale,
and a continuous mapping argument is used to complete the proof.
\end{abstract}



\end{frontmatter}

\section{Introduction}
\label{secintroduction}

Sampling probability distributions $\pi^N$ in $\RR^N$ for $N$ large
is of interest in numerous applications arising in
applied probability and statistics. The Markov chain Monte Carlo
(MCMC) methodology~\cite{CaseRobe04} provides a framework for many algorithms
which affect this sampling. It is hence of interest to quantify
the computational cost of MCMC methods as a function of
dimension $N$. This paper is part of a research program designed to
develop the analysis of MCMC in high dimensions so that it may be usefully
applied to understand target measures which arise in applications.
The simplest class of target measures for
which analysis can be carried out are perhaps target
distributions $\pi^N$ of the form
\begin{equation}
\label{eqtarget1}
\frac{d\pi^{N}}{d \lambda^N}(x)=\prod_{i=1}^N f(x_i).
\end{equation}
Here $\lambda^N(dx)$ is the $N$-dimensional Lebesgue measure,
and $f(x)$ is a one-dimensional probability density function.
Thus $\pi^N$ has the form of an i.i.d. product.
Using understanding gained in this situation, we will
develop an analysis, that is, relevant to an important class of
nonproduct measures which arise in a range of applications.

We start by describing the MCMC methods which are studied in this paper.
Consider a $\pi^N$-invariant metropolis Hastings--Markov chain
$\{ x^{k,N} \}_{k \geq1}.$ From the current state $x$,
we propose $y$ drawn from the kernel $q(x,y);$ this is
then accepted with probability
\[
\alpha(x, y) = 1 \wedge{\pi^N(y) q(y,x) \over\pi^N(x) q(x,y)}.
\]
Two widely used proposals are the random walk proposal
(obtained from the discrete approximation of Brownian motion)
\begin{equation} \label{eqnRWMdef}
y = x + \sqrt{2\delta} Z^N,\qquad Z^N \sim\Normal(0, \mathrm{I}_N),
\end{equation}
and the Langevin proposal
(obtained from the time discretization of the Langevin diffusion)
\begin{equation} \label{eqnMALAalgdef}
y = x + \delta\nabla\log\pi^N(x) + \sqrt{2 \delta} Z^N,\qquad
Z^N \sim\Normal(0, \mathrm{I}_N) .
\end{equation}
Here $2\delta$ is the proposal variance,
a parameter quantifying the size of the discrete time increment; we
will consider ``local proposals'' for which $\delta$ is small.
The Markov chain corresponding to proposal \eqref{eqnRWMdef} is the
Random Walk Metropolis (RWM) algorithm~\cite{Metretal53},
and the Markov transition rule constructed from
the proposal \eqref{eqnMALAalgdef} is known as the Metropolis Adjusted Langevin
Algorithm (MALA)~\cite{CaseRobe04}. This paper is aimed
at analyzing the computational complexity of
the MALA algorithm in high dimensions.

A fruitful way to quantify the computational cost of these Markov chains
which proceed via local\vadjust{\goodbreak} proposals is to determine
the ``optimal'' size of increment $\delta$ as a function of dimension $N$
(the precise notion of optimality is discussed below).
A~simple heuristic suggests the existence of such an ``optimal scale''
for $\delta$: smaller values of
the proposal variance lead to high acceptance rates, but the chain does
not move much even when accepted,
and therefore may not be efficient. Larger values of the proposal variance
lead to larger moves, but then the
acceptance probability is tiny. The optimal scale for
the proposal variance strikes a balance between making large moves and
still having
a reasonable acceptance probability.
In order to quantify this idea it is useful to define
a continuous interpolant of the Markov chain as follows:
\begin{eqnarray} \label{eqnMCMCe}
z^N(t)= \biggl(\frac{t}{\Delta t} - k\biggr) x^{k+1,N}
+ \biggl(k+1 - \frac{t}{\Delta t}\biggr) x^{k,N}
\nonumber
\\[-8pt]
\\[-8pt]
\eqntext{\mbox{for } k \Delta t \leq t < (k+1) \Delta t.}
\end{eqnarray}
We choose the proposal variance to satisfy $\delta= \ell\Delta t$,
with $\Delta t = N^{-\gamma}$ setting the scale in terms of
dimension and the parameter
$\ell$ a ``tuning'' parameter which
is independent of the dimension $N$.
Key questions, then, concern the choice of $\gamma$ and
$\ell$. If $z^N$ converges weakly to a suitable stationary
diffusion process, then it is natural to deduce that
the number of Markov chain steps required in stationarity
is inversely proportional to the proposal variance,
and hence to $\Delta t,$ and so grows like~$N^{\gamma}$. The parametric dependence of the
limiting diffusion process then provides a selection
mechanism for $\ell$.
A~research program along these lines was initiated
by Roberts and coworkers in the pair of papers
\cite{Robeetal97,RobeRose98}.
These papers concerned the RWM and MALA algorithms,
respectively, when applied to the target \eqref{eqtarget1}.
In both cases it was shown that the projection of $z^N$ into any single
fixed coordinate direction $x_i$ converges weakly
in $C([0,T];\RR)$ to $z$, the scalar diffusion process
\begin{equation}
\frac{dz}{dt} = h(\ell) [\log f(z)]'+\sqrt{2 h(\ell) }\,\frac
{dW}{dt}\label{eqsde}
\end{equation}
for $h(\ell)>0$, a constant determined by the parameter
$\ell$ from the proposal variance.
For RWM the scaling of the proposal variance
to achieve this limit is determined by the
choice $\gamma=1$~\cite{Robeetal97}, while
for MALA $\gamma=\frac13$~\cite{RobeRose98}.
The analysis shows that the number of steps required
to sample the target measure
grows as $\OO(N)$ for RWM, but only as $\OO(N^{1/3})$ for MALA.
This quantifies the efficiency gained by
use of MALA over RWM, and in particular from employing local moves
informed by the gradient of the logarithm of the target density.
A second important feature of the analysis is that it suggests
that the optimal choice of $\ell$ is that which
maximizes $h(\ell)$. This value of $\ell$ leads, in both
cases, to a universal [independent of $f(\cdot)$]
optimal average acceptance probability
(to three significant figures)
of $0.234$ for RWM and $0.574$ for MALA.

These theoretical analyses have had a huge practical impact as the
optimal acceptance probabilities send
a concrete message to practitioners: one should
``tune'' the proposal variance of the RWM and MALA algorithms
so as to have acceptance probabilities of $0.234$ and $0.574$,
respectively.
However, practitioners use these tuning criteria far outside
the class of target distributions given by \eqref{eqtarget1}.
It is natural to ask whether they are wise to do so.
Extensive simulations (see~\cite{RobeRose01,Sheretal10}) show that these
optimality results also hold for more complex target distributions.
Furthermore, a range of subsequent theoretical analyses
confirmed that the optimal scaling ideas do indeed
extend beyond \eqref{eqtarget1}; these papers studied
slightly more complicated models, such as products of one-dimensional
distributions with different variances and elliptically
symmetric distributions
\cite{Beda07,Beda09,Breyetal04,christensen2005scaling}.
However, the diffusion limits obtained remain essentially one dimensional
in all of these extensions.\setcounter{footnote}{3}\footnote{The paper~\cite{BreyRobe00}
contains an infinite dimensional diffusion limit, but we have been unable
to employ the techniques of that paper.} In this paper we study considerably
more complex target distributions which are not of the product form,
and the limiting diffusion takes values in an infinite dimensional space.

Our perspective on these problems is motivated by applications
such as Bayesian nonparametric statistics, for example, in
application to inverse problems~\cite{Stua10}, and
the theory of conditioned diffusions~\cite{HairStuaVoss10}.
In both these areas the target measure of interest,
$\pi$, is on an infinite dimensional real separable Hilbert
space $\h$ and, for Gaussian priors (inverse problems)
or additive noise (diffusions)
is absolutely continuous with respect to a
Gaussian measure $\pi_0$ on $\h$
with mean zero and covariance operator $\C$.
This framework for the analysis of MCMC in high dimensions
was first studied in the papers~\cite{Besketal08,BRS09,BeskStua07}.
The Radon--Nikodym derivative defining the target measure is
assumed to have the form
\begin{equation}
\frac{d\pi}{d\pi_0}(x) = M_{\Psi} \exp( -\Psi(x) )
\label{eqntargmeas}
\end{equation}
for a real-valued functional
$\Psi\dvtx\h^s \mapsto\mathbb{R}$
defined on a subspace $\h^s \subset\h$ that
contains the support of the reference measure $\pi_0$;
here $M_{\Psi}$ is a normalizing constant.
We are interested in studying MCMC methods applied to
finite dimensional approximations of this measure
found by projecting onto the first $N$ eigenfunctions
of the covariance operator $\C$ of the Gaussian
reference measure $\pi_0$.

It is proved in
\cite{DaprZaby92,Hairetal05,HairStuaVoss07} that
the measure $\pi$ is invariant for $\mathcal{H}$-valued SDEs
(or stochastic PDEs--SPDEs) with the form
\begin{equation} \label{eqnspde}
\frac{dz}{dt} = - h(\ell) \bigl(z + \C\nabla\Psi(z) \bigr)+ \sqrt{2
h(\ell)}\, \frac{dW}{dt},\qquad
z(0)=z^0,
\end{equation}
where $W$ is a Brownian motion (see~\cite{DaprZaby92}) in
$\h$ with covariance operator $\C$.
In~\cite{MatPillStu09} the RWM\vadjust{\goodbreak} algorithm is studied
when applied to a sequence of finite dimensional approximations
of $\pi$ as in \eqref{eqntargmeas}.
The continuous time interpolant
of the Markov chain $z^N$ given by \eqref{eqnMCMCe} is
shown to converge
weakly to $z$ solving \eqref{eqnspde} in $C([0,T];\h^s)$.
Furthermore, as for the i.i.d. target measure, the scaling
of the proposal variance which achieves this scaling
limit is inversely proportional to $N$ (i.e., corresponds to the exponent
$\gamma=1$), and the speed
of the limiting diffusion process is maximized at the
same universal acceptance probability of $0.234$ that
was found in the i.i.d. case. Thus, remarkably, the i.i.d.
case has been of fundamental importance in understanding
MCMC methods applied to complex infinite dimensional
probability measures arising in practice.
The paper~\cite{MatPillStu09} developed an approach for
deriving diffusion limits for such algorithms, using
ideas from numerical analysis. We can build on these techniques to
derive scaling limits for
a wide range of Metropolis--Hastings algorithms
with local proposals.

The purpose of this article is to develop the techniques
in the context of the MALA algorithm. To the best of our
knowledge, the only paper to consider the optimal
scaling for the MALA algorithm for nonproduct targets
is~\cite{Breyetal04}, in the context of nonlinear regression.
In~\cite{Breyetal04} the target measure has a structure similar
to that of the mean field models studied in statistical
mechanics and hence behaves
asymptotically like a product measure when the
dimension goes to infinity. Thus the diffusion limit obtained
in~\cite{Breyetal04} is finite dimensional.

The main contribution of our work is the proof of a diffusion limit
for the output of the MALA algorithm, suitably
interpolated, to the SPDE
\eqref{eqnspde}, when applied to
$N$-dimensional approximations
of the target measures \eqref{eqntargmeas}
with proposal variance inversely
proportional to $N^{1/3}$. Moreover we show that the
speed $h(\ell)$ of the limiting diffusion is maximized
for an average acceptance probability of $0.574$, just
as in the i.i.d. product scenario~\cite{RobeRose98}.
Thus in this regard,
our work is the first extension of the remarkable
results in~\cite{RobeRose98} for the Langevin algorithm
to target measures which are not of product form. This adds
theoretical weight to the results observed in computational experiments
which demonstrate the robustness of the optimality
criteria developed in~\cite{Robeetal97,RobeRose98}.
In particular, the paper~\cite{Besketal08} shows
numerical results indicating the need to scale
time-step as a function of dimension to obtain
${\mathcal O}(1)$ acceptance probabilities.

In Section~\ref{secmainthm} we state the main theorem of the
paper, having defined precisely the setting in which it holds.
Section~\ref{secproof} contains the proof of the main theorem,
postponing the proof of a number of key technical estimates to
Section~\ref{seckeyestimates}. In Section~\ref{secconclusion}
we conclude by summarizing and providing the outlook for further
research in this area.

\section{Main theorem}
\label{secmainthm}

This section is devoted to stating the main theorem of the
article. However, the setting is complex, and we develop it in a
step-by-step fashion, before the theorem statement.
In Section~\ref{secprior} we introduce the form of the
reference, or prior, Gaussian measure $\pi_0$, followed
in Section~\ref{sseccom} by the change of measure
which induces a genuinely nonproduct structure. In
Section~\ref{ssecapprox} we describe finite dimensional
approximation of the measure, enabling us to define
application of a variant MALA-type algorithm
in Section~\ref{secalgorithm}. We then discuss
in Section~\ref{secoptimalscale}
how the choice of scaling used in the theorem emerges from
study of the acceptance probabilities. Finally, in Section
\ref{ssecmt}, we state the main theorem.

Throughout the paper we use the
following notation in order to compare
sequences and to denote conditional expectations:
\begin{itemize}
\item Two sequences $\{\alpha_n\}$ and $\{\beta_n\}$ satisfy $\alpha_n
\lesssim\beta_n$
if there exists a constant $K>0$ satisfying $\alpha_n \leq K \beta_n$
for all $n \geq0$.
The notation $\alpha_n \asymp\beta_n$ means that $\alpha_n \lesssim
\beta_n$ and $\beta_n \lesssim\alpha_n$.
\item
Two sequences of real functions $\{f_n\}$ and $\{g_n\}$ defined on the
same set $D$
satisfy $f_n \lesssim g_n$ if there exists a constant $K>0$ satisfying
$f_n(x) \leq K g_n(x)$ for all $n \geq0$
and all $x \in D$.
The notation $f_n \asymp g_n$ means that $f_n \lesssim g_n$ and $g_n
\lesssim f_n$.
\item
The notation $\EE_x [ f(x,\xi) ]$ denotes
expectation with respect to $\xi$ with
the variable~$x$ fixed.
\end{itemize}

\subsection{Gaussian reference measure}
\label{secprior}
Let $\h$ be a separable Hilbert space of real valued functions
with scalar product denoted by $\langle\cdot, \cdot \rangle$
and associated norm $\|x\|^2 = \langle x,x \rangle$.
Consider a Gaussian probability measure $\pi_0$
on $(\h, \| \cdot\| )$ with covariance
operator~$\C$. The general theory of Gaussian
measures~\cite{DaprZaby92} ensures
that the operator~$\C$ is positive and trace class.
Let $\{\phi_j,\lambda^2_j\}_{j \geq1}$ be the eigenfunctions
and eigenvalues of the covariance operator $\C$:
\[
\C\phi_j = \lambda^2_j \phi_j,\qquad j \geq1.
\]
We assume a normalization under which
the family $\{\phi_j\}_{j \geq1}$
forms a complete orthonormal basis in
the Hilbert space $\h$, which
we refer to us as the Karhunen--Lo\`eve basis.
Any function $x \in\h$ can be represented in this
basis via the expansion
\begin{equation}
x= \sum_{j=1}^{\infty} x_j \phi_j,\qquad x_j \eqdef\langle x,\phi_j
\rangle.
\label{eqneigenexp}
\end{equation}
Throughout this paper we will often identify the function $x$ with
its coordinates $\{x_j\}_{j=1}^{\infty} \in\ell^2$ in this
eigenbasis, moving freely between the\vspace*{1pt} two representations.
The Karhunen--Lo\`eve expansion (see~\cite{DaprZaby92},
Section \textit{White noise expansions}),
refers to the fact that a realization $x$ from the Gaussian measure
$\pi
_0$ can
be expressed by allowing the coordinates $\{x_j\}_{j \geq1}$ in
\eqref{eqneigenexp} to be independent random
variables distributed as $x_j \sim\Normal(0,\lambda_j^2)$. Thus, in
the coordinates $\{x_j\}_{j \geq1}$,
the Gaussian reference measure $\pi_0$
has a product structure.

For every $x \in\h$ we have representation
\eqref{eqneigenexp}.
Using this expansion, we define Sobolev-like spaces $\h^r, r \in\RR$,
with the inner-products and norms defined by
%
%
\begin{equation}\label{eqnSob}
\langle x,y \rangle_r \eqdef\sum_{j=1}^\infty j^{2r}x_jy_j,\qquad
\|x\|^2_r \eqdef\sum_{j=1}^\infty j^{2r} x_j^{2}.
\end{equation}
Notice that $\h^0 = \h$ and
$\h^r \subset\h\subset\h^{-r}$ for any $r >0$.
The Hilbert--Schmidt norm $\|\cdot\|_\C$ associated to the covariance
operator $\C$
is defined as
\[
\|x\|^2_\C= \sum_j \lambda_j^{-2} x_j^2.
\]
For $x,y \in\h^r$, the outer product operator in $\h^r$
is the operator $x \otimes_{\h^r} y\dvtx\h^r \to\h^r$
defined by $(x \otimes_{\h^r} y) z \eqdef\langle y, z \rangle_r x$
for every $z \in\h^r$.
For $r \in\RR$, let $B_r \dvtx\h\mapsto\h$ denote the operator
which is
diagonal in the basis $\{\phi_j\}_{j \geq1}$ with diagonal entries
$j^{2r}$. The operator $B_r$ satisfies
$B_r \phi_j = j^{2r} \phi_j$
so that $B^{1/2}_r \phi_j = j^r \phi_j$.
The operator $B_r$
lets us alternate between the Hilbert space $\h$ and the
Sobolev spaces $\h^r$ via the identities
$\langle x,y \rangle_r = \langle B^{1/2}_r x,B^{1/2}_r y
\rangle$.
Since $\|B_r^{-1/2} \phi_k\|_r = \|\phi_k\|=1$,
we deduce that $\{B^{-1/2}_r \phi_k \}_{k \geq0}$ forms an
orthonormal basis for $\h^r$.
For a positive, self-adjoint operator $D \dvtx\h\mapsto\h$,
we define its trace in $\h^r$ by
\begin{equation}
\label{eqntrace}
\tr_{\h^r}(D)
\eqdef\sum_{j=1}^\infty\langle(B_r^{-{1}/{2}} \phi_j), D
(B_r^{-{1}/{2}} \phi_j) \rangle_r.
\end{equation}
Since $\tr_{\h^r}(D)$ does not depend on the orthonormal basis, 
the operator $D$ is said to be trace class in $\h^r$ if $\tr_{\h^r}(D)
< \infty$ for
some, and hence any, orthonormal basis of $\h^r$.
Let us define the operator $\C_r \eqdef B^{1/2}_r \C B^{1/2}_r$.
Notice that $ \tr_{\h^r}(\C_r)=\sum_{j=1}^\infty\lambda^2_j j^{2r}$.
In~\cite{MatPillStu09} it is shown that under the condition
\begin{equation}
\label{eqfinitetrace}
\tr_{\h^r}(\C_r) < \infty,
\end{equation}
the support of $\pi_0$ is included in $\h^r$ in the sense that
$\pi_0$-almost every function $x \in\h$ belongs to $\h^r$.
Furthermore, the induced distribution of $\pi_0$ on $\h^r$ is
identical to that of a centered Gaussian measure
on $\h^r$ with covariance operator $\C_r$.
For example, if $\xi\dist\pi_0$,
then $\EE[ \langle\xi,u \rangle_r \langle\xi,v \rangle_r ] =
\langle u,\C_r v \rangle_r$
for any functions $u,v \in\h^r$.
Thus in what follows, we alternate between
the Gaussian measures $\Normal(0,\C)$ on $\h$
and $\Normal(0,\C_r)$ on $\h^r$, for those $r$ for
which \eqref{eqfinitetrace} holds.

\subsection{Change of measure}
\label{sseccom}

Our goal is to sample from a measure $\pi$
defined through the change of probability formula
\eqref{eqntargmeas}. As described in Section~\ref{secprior},
the condition $\tr_{\h^r}(\C_r) < \infty$ implies that
the measure $\pi_0$ has full support on $\h^r$,
that is, $\pi_0(\h^r)=1$. Consequently,
if $\tr_{\h^r}(\C_r) < \infty$, the functional
$\Psi(\cdot)$ needs only to be defined on $\h^r$
in order for the change of probability formula
\eqref{eqntargmeas} to be valid. In this section,
we give assumptions on the decay of the
eigenvalues of the covariance operator $\C$ of $\pi_0$
that ensure the existence of a real
number $s>0$ such that $\pi_0$ has full support
on $\h^s$. The functional $\Psi(\cdot)$ is assumed
to be defined on $\h^s$, and we impose regularity
assumptions on $\Psi(\cdot)$ that ensure that the
probability distribution $\pi$ is not too different
from $\pi_0$, when projected into directions associated
with $\varphi_j$ for $j$ large.
For each $x \in\h^s$ the derivative $\nabla\Psi(x)$
is an element of the dual $(\h^s)^*$ of $\h^s$,
comprising linear functionals on $\h^s$.
However, we may identify $(\h^s)^*$ with $\h^{-s}$
and view $\nabla\Psi(x)$
as an element of $\h^{-s}$ for each $x \in\h^s$. With this identification,
the following identity holds:
\[
\| \nabla\Psi(x)\|_{\mathcal{L}(\h^s,\RR)} = \| \nabla\Psi(x) \|_{-s},
\]
and the second derivative $\partial^2 \Psi(x)$ can
be identified as an element
of $\mathcal{L}(\h^s, \h^{-s})$.
To avoid technicalities we assume that $\Psi(\cdot)$ is quadratically bounded,
with the first derivative linearly bounded, and the second derivative globally
bounded. Weaker assumptions could be dealt with by use of stopping time
arguments.
%
\begin{assumption} \label{ass1}
The covariance operator $\C$ and functional $\Psi$ satisfy the following:
\begin{longlist}[(1)]
\item[(1)] \textit{Decay of Eigenvalues $\lambda_j^2$ of $\C$}:
there is an exponent $\kappa> \frac{1}{2}$ such that
\begin{equation}
\label{eqeigenvaluedecay}
\lambda_j \asymp j^{-\kappa}.
\end{equation}
\item[(2)] \textit{Assumptions on} $\Psi$:
There exist constants $M_i \in\RR, i \leq4$, and $s \in
[0, \kappa- 1/2)$ such that for all $x \in\h^s$ the functional
$\Psi\dvtx\h^s \to\RR$ satisfies
\begin{eqnarray}
M_1 &\leq&\Psi(x) \leq M_2 (1 + \|x\|_s^2), \label{eqnpsi1}\\
\| \nabla\Psi(x)\|_{-s} &\leq &M_3 (1 + \|x\|_s), \label{eqnpsi2}\\
\opnm{\partial^2 \Psi(x)}{s}{-s} & \leq& M_4. \label{eqnpsi3}
\end{eqnarray}
\end{longlist}
\end{assumption}

\begin{remark}
The condition $\kappa>\frac{1}{2}$ ensures that the covariance
operator $\C$
is trace class in $\h$.
In fact, equation \eqref{eqfinitetrace} shows that $\C_r$
is trace-class in $\h^r$ for any $r < \kappa- \frac{1}{2}$.
It follows that $\pi_0$ has full measure in $\h^r$
for any $r \in[0, \kappa- 1/2)$. In particular $\pi_0$ has full
support on $\h^s$.
\end{remark}
\begin{remark}
The functional $\Psi(x) = \frac{1}{2}\|x\|_s^2$
satisfies Assumption~\ref{ass1}.
It is defined on $\h^s$ and its derivative at $x \in\h^s$
is given by $\nabla\Psi(x) = \sum_{j \geq0} j^{2s} x_j \phi_j \in
\h
^{-s}$ with
$\|\nabla\Psi(x)\|_{-s} = \|x\|_s$.
The second derivative $\partial^2 \Psi(x) \in\mathcal{L}(\h^s, \h^{-s})$
is the linear operator that maps $u \in\h^s$ to $\sum_{j \geq0}
j^{2s} \langle u,\phi_j \rangle \phi_j \in\h^s$:
its norm satisfies $\| \partial^2 \Psi(x) \|_{\mathcal{L}(\h^s, \h
^{-s})} = 1$ for any $x \in\h^s$.
\end{remark}

Since the eigenvalues $\lambda_j^2$ of $\C$ decrease as $\lambda_j
\asymp j^{-\kappa}$,
the operator $\C$ has a smoothing effect: $\C^{\alpha} h$ gains $2
\alpha\kappa$ orders of regularity
in the sense that the $\h^{\beta}$-norm of $\C^{\alpha}h$ is
controlled by
the $\h^{\beta-2 \alpha\kappa}$-norm of $h \in\h$.
Indeed, under Assumption~\ref{ass1}, the following estimates holds:
\begin{equation}
\label{eqgainregularity}
\|h\|_{\C} \asymp\|h\|_{\kappa}
\quad\mbox{and}\quad
\| \C^{\alpha}h \|_{\beta} \asymp\|h\|_{\beta- 2 \alpha\kappa}.
\end{equation}
The proof follows the methodology used to prove Lemma $3.3$ of
\cite{MatPillStu09}. The reader is referred to this text for more details.

\subsection{Finite dimensional approximation}
\label{ssecapprox}
We are interested in finite dimensional approximations
of the probability distribution $\pi$.
To this end, we introduce the vector space spanned by
the first $N$ eigenfunctions of the covariance operator,
\[
X^N \eqdef{\mathrm{span}}\{ \phi_1, \phi_2, \ldots, \phi_N \}.
\]
Notice that $X^N \subset\h^r$ for any $r \in[0; +\infty)$. In particular,
$X^N$ is a subspace of~$\h^s$.
Next, we define $N$-dimensional approximations
of the functional $\Psi(\cdot)$
and of the reference measure $\pi_0$. To this end,
we introduce the orthogonal projection on $X^N$
denoted by $P^N \dvtx\h^s \mapsto X^N \subset\h^s$.
The functional $\Psi(\cdot)$ is approximated by the
functional $\Psi^N\dvtx X^N \mapsto\RR$ defined by
\begin{equation}
\label{eqpsiN}
\Psi^N \eqdef\Psi\circ P^N.
\end{equation}
The approximation $\pi_0^N$ of the reference measure $\pi_0$ is the Gaussian
measure on $X^N$ given by the law of the random variable
\[
\pi_0^N \dist\sum_{j=1}^N \lambda_j \xi_j \phi_j
= (\C^N)^{1/2} \xi^N,
\]
where $\xi_j$ are i.i.d. standard Gaussian random variables, $\xi^N =
\sum_{j=1}^N \xi_j \phi_j$
and $\C^N = P^N \circ\C\circ P^N$. Consequently we have
$\pi_0^N = \Normal(0, \C^N)$.
Finally, one can define the approximation $\pi^N$
of $\pi$ by the change of probability formula
\begin{equation} \label{eqtarget3}
\frac{d\pi^{N}}{d\pi^N_0}(x) = M_{\Psi^N} \exp(-\Psi^N(x)),
\end{equation}
where $M_{\Psi^N}$ is a normalization constant. Notice that
the probability distribution $\pi^N$ is supported on $X^N$
and has Lebesgue density\footnote{For ease of notation we do not distinguish
between a measure and its density, nor do we
distinguish between the representation of the
measure in $X^N$ or in coordinates in $\RR^N$.}
on $X^N$ equal to
\begin{equation}\label{eqnpitrunc}
\pi^N(x) \propto
\exp\bigl( -\tfrac{1}{2} \|x\|^2_{\C^N} -\Psi^N(x) \bigr).
\end{equation}
In formula \eqref{eqnpitrunc},
the Hilbert--Schmidt norm\vspace*{1pt} $\| \cdot\|_{\C^N}$ on $X^N$ is given by the
scalar product
$\langle u,v \rangle_{\C^N} = \langle u, (\C^N)^{-1} v \rangle$ for
all $u,v \in X^N$.
The operator $\C^{N}$ is invertible on $X^N$ because the
eigenvalues of $\C$ are assumed to be strictly positive.
The quantity $\C^N \nabla\log\pi^N(x)$ is repeatedly used in the
text and,
in particular, appears in the function $\mu^N(x)$ given by
\begin{equation} \label{eqnmun}
\mu^N(x) = -\bigl( P^N x + \C^N \nabla\Psi^N(x) \bigr)
\end{equation}
which, up to an additive constant, is $\C^N \nabla\log\pi^N(x).$
This function is the drift of an ergodic Langevin diffusion
that leaves $\pi^N$ invariants. Similarly, one defines the
function $\mu\dvtx\h^s \to\h^s$ given by
\begin{equation}
\label{eqnmu}
\mu(x) = - \bigl( x + \C\nabla\Psi(x) \bigr)
\end{equation}
which can informally be seen as
$\C\nabla\log\pi(x)$, up to an additive constant.
In the sequel, Lemma~\ref{lemmuNmu} shows that, for $\pi_0$-almost every
function $x \in\h$, we have $\lim_{N \to\infty} \mu^N(x) = \mu(x)$.
This quantifies the manner in which $\mu^N(\cdot)$ is an
approximation of
$\mu(\cdot)$.

The next lemma gathers
various regularity estimates on the functional $\Psi(\cdot)$
and $\Psi^N( \cdot)$ that are repeatedly used in the sequel.
These are simple consequences of Assumption~\ref{ass1}, and proofs
can be found in~\cite{MatPillStu09}.
\begin{lemma}[(Properties of $\Psi$)] \label{lemregularity}
Let the functional $\Psi(\cdot)$ satisfy Assumption~\ref{ass1} and
consider the functional $\Psi^N(\cdot)$ defined by equation \eqref{eqpsiN}.
The following estimates hold:
\begin{longlist}[(1)]
\item[(1)]
The functionals $\Psi^N\dvtx\h^s \to\RR$ satisfy the same conditions
imposed on $\Psi$ given by equations \eqref{eqnpsi1},
\eqref{eqnpsi2} and \eqref{eqnpsi3}
with constants that can be chosen independent of $N$.
\item[(2)] The function $\C\nabla\Psi\dvtx\h^s \to\h^s$ is globally
Lipschitz on $\h^s$:
there exists a constant $M_5 > 0$ such that
\[
\| \C\nabla\Psi(x) - \C\nabla\Psi(y) \|_s \leq M_5 \| x-y \|_s\qquad
\forall x,y \in\h^s.
\]
Moreover, the functions $\C^N \nabla\Psi^N\dvtx\h^s \to\h^s$
also satisfy
this estimate with a constant that can be chosen independently of $N$.
\item[(3)]
The functional $\Psi(\cdot)\dvtx\h^s \to\RR$ satisfies a second order
Taylor formula.\footnote{We extend $\langle\cdot,\cdot\rangle$ from an
inner-product on $\h$ to the dual pairing between
$\h^{-s}$ and $\h^s$.} There exists a constant $M_6 > 0$ such that
\begin{equation} \label{eqn2nd-Taylor}\quad
\Psi(y) - \bigl( \Psi(x) + \langle\nabla\Psi(x), y-x \rangle \bigr)
\leq M_6 \|x - y \|_s^2\qquad
\forall x,y \in\h^s.
\end{equation}
Moreover, the functionals $\Psi^N(\cdot)$ also satisfy
this estimates with a constant that can be chosen independently of $N$.
\end{longlist}
\end{lemma}

\begin{remark}
\label{remts}
Regularity Lemma~\ref{lemregularity} shows, in particular, that the
function $\mu\dvtx\h^s \to\h^s$ defined by \eqref{eqnmu}
is globally Lipschitz on $\h^s$.
Similarly, it follows that $\C^N \nabla\Psi^N\dvtx\h^s \to\h^s$
and $\mu^N\dvtx\h^s \to\h^s$ given by \eqref{eqnmun}
are globally Lipschitz with Lipschitz constants that
can be chosen uniformly in $N$.
\end{remark}

\subsection{The algorithm}\label{secalgorithm}
The MALA algorithm is defined in this section.
This method is motivated by the fact that the probability
measure $\pi^N$ defined by equation \eqref{eqtarget3}
is invariant with respect to the Langevin diffusion process
\begin{equation} \label{eqnspdeN}
\frac{dz}{dt} = \mu^N(z)+ \sqrt{2} \,\frac{dW^N}{dt},
\end{equation}
where $W^N$ is a Brownian motion in
$\h$ with covariance operator $\C^N$. The drift
function $\mu^N\dvtx\h^s \to\h^s$ is the gradient of
the log-density of $\pi^N$, as described by equation \eqref{eqnmun}.
The idea of the MALA algorithm is to make a proposal based
on Euler--Maruyama discretization of the diffusion \eqref{eqnspdeN}.
To this end we consider, from state $x \in X^N$, proposals $y \in X^N$
given by
\begin{equation} \label{eqnproposal}
y - x = \delta\mu^N(x) + \sqrt{2 \delta} (\C^N)^{{1}/{2}}\xi^N\qquad
\mbox{where }
\delta= \ell N^{-{1}/{3}}
\end{equation}
with $\xi^N= \sum_{i=1}^N \xi_i \phi_i$ and $\xi_i \dist\Normal(0,1)$.
Notice that $(\C^N)^{{1}/{2}}\xi^N \dist\Normal(0, \C^N)$.
The quantity $\delta$ is the time-step in an Euler--Maruyama
discretization of \eqref{eqnspdeN}.
We introduce a related parameter
\[
\Delta t:=\ell^{-1}\delta=N^{-1/3}
\]
which will be the natural time-step for the
limiting diffusion process derived from the proposal
above, after inclusion of an accept--reject mechanism.
The scaling of $\Delta t$, and hence $\delta,$
with $N$ will ensure
that the average acceptance probability is of order $1$
as $N$ grows. This is discussed in more detail in Section
\ref{secoptimalscale}. The quantity $\ell> 0$ is
a fixed parameter which can be
chosen to maximize the speed of the limiting diffusion process;
see the discussion in the \hyperref[secintroduction]{Introduction} and after the
Main Theorem below.

We will study the Markov chain $x^N = \{x^{k,N}\}_{k \geq0}$
resulting from Metropolizing this proposal when it is
started at stationarity:
the initial position $x^{0,N}$ is distributed as $\pi^N$ and
thus lies in $X^N$. Therefore, the Markov chain evolves
in $X^N$; as a consequence, only the first $N$ components of
an expansion in the eigenbasis of $\C$
are nonzero, and the algorithm can be implemented in $\RR^N$. However
the analysis is cleaner
when written in $X^N \subset\h^s$.
The acceptance probability only depends on
the first $N$ coordinates of $x$ and $y$ and has the form
\begin{equation} \label{eqnaccprob}
\alpha^N(x,\xi^N) = 1 \wedge\frac{\pi^N(y) T^N(y,x)}{\pi^N(x)
T^N(x,y)} = 1 \wedge e^{ Q^N(x,\xi^N) },
\end{equation}
where the proposal $y$ is given by equation \eqref{eqnproposal}.
The function $T^N(\cdot, \cdot)$ is the density of the Langevin
proposals \eqref{eqnproposal}
and is given by
\[
T^N(x,y) \propto\exp\biggl\{ -\frac{1}{4 \delta} \| y - x - \delta\mu
^N(x) \|^2_{\C^N} \biggr\}.
\]
The local mean acceptance probability $\alpha^N(x)$ is
defined by
\begin{equation} \label{elocaccept}
\alpha^N(x) = \EE_x[ \alpha^N(x,\xi^N) ].
\end{equation}
It is the expected acceptance probability when the algorithm stands
at $x \in\h$.
The Markov chain $x^N = \{x^{k,N}\}_{k \geq0}$ can also be expressed as
\begin{equation} \label{eqnMALAalgorithm}
\cases{
y^{k,N} = x^{k,N} + \delta\mu^N(x^{k,N}) +
\sqrt{2 \delta} (\C^N)^{{1}/{2}} \xi^{k,N} ,\vspace*{2pt}\cr
x^{k+1,N} = \gamma^{k,N} y^{k,N} + (1- \gamma^{k,N}) x^{k,N}, }
\end{equation}
where $\xi^{k,N}$ are i.i.d. samples distributed as $\xi^N$, and
$\gamma^{k,N} = \gamma^N(x^{k,N},\xi^{k,N})$ creates a
Bernoulli random sequence
with $k^{th}$ success probability $\alpha^N(x^{k,N},\break\xi^{k,N})$.
We may view the Bernoulli random variable as
$\gamma^{k,N} = 1_{\{U^k < \alpha^N(x^{k,N},\xi^{k,N})\}}$ where
$U^k \dist\operatorname{Uniform}(0,1)$ is independent from $x^{k,N}$
and $\xi^{k,N}$.
The quantity $Q^N$ defined in equation \eqref{eqnaccprob} may be
expressed as
\begin{eqnarray} \label{eqnQN}
\qquad Q^N(x, \xi^N)
&=& -\frac{1}{2} ( \|y\|_{\C^N}^2 - \|x\|_{\C^N}^2 ) -\bigl( \Psi
^N(y)-\Psi^N(x) \bigr)
\nonumber
\\[-8pt]
\\[-8pt]
\nonumber
&&{} -\frac{1}{4 \delta} \{
\| x - y - \delta\mu^N(y) \|^2_{\C^N}
- \| y - x -\delta\mu^N(x) \|^2_{\C^N} \}.
\end{eqnarray}
As will be seen in the next section, a key idea behind our
diffusion limit is that, for large $N$, the quantity $Q^N(x,\xi^N)$
behaves like
a Gaussian random variable
independent from the current position $x$.

In summary, the Markov chain that we have described in $\h^s$ is,
when projected onto $X^N$, equivalent to a standard
MALA algorithm on $\RR^N$ for the Lebesgue density \eqref{eqnpitrunc}.
Recall that the target measure $\pi$ in \eqref{eqntargmeas} is the
invariant measure of the SPDE \eqref{eqnspde}. Our goal is to obtain an
invariance principle for the continuous interpolant \eqref{eqnMCMCe}
of the Markov chain $x^N = \{x^{k,N}\}_{k \geq0}$ started in stationarity,\vspace*{1pt}
that is, to show weak convergence in $C([0,T]; \h^s)$ of
$z^{N}(t)$ to the solution $z(t)$ of the SPDE \eqref{eqnspde},
as the dimension $N \rightarrow\infty$.

\subsection{\texorpdfstring{Optimal scale $\gamma=\frac13$}{Optimal scale gamma=1/3}}
\label{secoptimalscale}
In this section, we informally describe why the optimal scale
for the MALA proposals \eqref{eqnproposal}
is given by the exponent $\gamma= \frac{1}{3}$.
For product-form target probability
described by equation \eqref{eqtarget1}, the optimality of the exponent
$\gamma= \frac{1}{3}$ was first obtained in~\cite{RobeRose98}.
For further discussion, see also~\cite{BRS09}. To keep the
exposition simple in this explanatory subsection,
we focus on the case $\Psi(\cdot) = 0$. The analysis is similar
with a nonvanishing function $\Psi(\cdot)$, because absolute
continuity ensures that the effect of $\Psi(\cdot)$ is
small compared to the dominant Gaussian effects described here.
Inclusion of nonvanishing $\Psi(\cdot)$
is carried out in Lemma~\ref{lemgaussianapprox}.

In the case
$\Psi(\cdot) = 0$, straightforward algebra shows that the acceptance
probability
$\alpha^N(x,\xi^N) = 1 \wedge e^{ Q^N(x,\xi^N) }$ satisfies
\[
Q^N(x,\xi^N)
= -\frac{\ell\Delta t}{4} ( \|y\|^2_{\C^N} - \|x\|^2_{\C^N}).
\]
For $\Psi(\cdot)=0$ and $x \in X^N$, the proposal $y$ is distributed as
$y = (1-\ell\Delta t)x + \sqrt{2 \ell\Delta t} (\C^N)^{1/2}
\xi^N$.
It follows that
\begin{eqnarray*}
\|y\|^2_{\C^N} - \|x\|^2_{\C^N}
&=&
-2 \ell\Delta t \bigl( \|x\|_{\C^N}^2 - \|(\C^N)^{1/2}\xi^N\|
^2_{\C^N}\bigr)
+ (\ell\Delta t)^2 \|x\|^2_{\C^N}\\
&&{}
+ 2\sqrt{2 \ell\Delta t} (1-\Delta t) \langle x, (\C^N)^{1/2}\xi ^N \rangle_{\C^N}.
\end{eqnarray*}
The details can be found in the proof of
Lemma~\ref{lemgaussianapprox}. Since the Markov chain $x^N=\{x^{k,N}\}
_{k \geq0}$
evolves in stationarity, for all $k \geq0$, we have $x^{k,N} \dist\pi
^N = \Normal(0,\C^N)$.
Therefore, with $x \dist\Normal(0,\C^N)$ and $\xi^N \dist\Normal
(0,\C^N)$,
the law of large numbers shows that both
$\|x\|_{\C^N}^2$ and $\|(\C^N)^{1/2}\xi^N\|^2_{\C^N}$
are of order $\OO(N)$,
while the central limit theorem shows that
$\langle x, (C^N)^{1/2}\xi^N \rangle_{\C^N} = \OO(N^{1/2})$
and $\|x\|_{\C^N}^2-\|(\C^N)^{1/2}\xi^N\|^2_{\C^N}=\OO
(N^{1/2})$.
For $\Delta t = \ell N^{-\gamma}$ and $\gamma< \frac{1}{3}$, it follows
\[
Q^N(x,\xi^N)
= -\frac{(\ell\Delta t)^3}{4}\|x\|^2_{\C^N} + \OO(N^{{1}/{2}-
{3 \gamma}/{2}})
\approx-\frac{\ell^3}{4} N^{1-3 \gamma},
\]
which shows that the acceptance probability is exponentially small of order
$\exp( -\frac{\ell^3}{4} N^{1-3 \gamma} )$. The same argument
shows that
for $\gamma> \frac{1}{3}$, we have $Q^N(x,\xi^N) \to0$, which shows
that the average
acceptance probability converges to $1$. For the critical exponent
$\gamma= \frac{1}{3}$,
the acceptance probability is of order $\OO(1)$.
In fact Lemma~\ref{lemgaussianapprox}
shows that for $\gamma= \frac13$, even when $\Psi(\cdot)$ is
nonzero, the following Gaussian approximation holds:
\[
Q^N(x,\xi^N) \approx\Normal\biggl(-\frac{\ell^3}{4}, \frac{\ell^3}{2}\biggr).
\]
This approximation is key to derivation of the diffusion limit.
In summary, choosing $\gamma> \frac{1}{3}$ leads to exponentially
small acceptance probabilities:
almost all the proposals are rejected so that the expected squared jumping
distance $\EE_{\pi^N}[\|x^{k+1,N}-x^{k,N}\|^2]$ converges
exponentially quickly
to $0$ as the dimension $N$ goes to infinity.
On the other hand, for any exponent $\gamma\geq\frac{1}{3}$,
the acceptance probabilities are bounded away from zero: the Markov chain
moves with jumps of size $\OO(N^{-{\gamma}/{2}})$, and the expected
squared jumping
distance is of order $\OO(N^{-\gamma})$. If we adopt the expected
squared jumping
distance as measure of efficiency,\vspace*{1pt} the optimal exponent is thus
given by $\gamma= \frac{1}{3}$. This viewpoint is analyzed
further in~\cite{BRS09}.

\subsection{Statement of main theorem}
\label{ssecmt}

The main result of this article describes the behavior of the MALA
algorithm for the
optimal scale $\gamma= \frac13$; the proposal variance is given by
$\delta= 2 \ell N^{-1/3}$.
In this case, Lemma~\ref{lemgaussianapprox} shows that the local mean
acceptance
probability $\alpha^N(x,\xi^N) = 1 \wedge e^{ Q^N(x,\xi^N) }$
satisfies $Q^N(x,\xi^N) \to Z_{\ell} \dist\Normal(-\frac{\ell^3}{4},
\frac{\ell^3}{2})$.
As a consequence, the asymptotic mean
acceptance probability of the MALA algorithm can
be explicitly computed as a function
of the parameter $\ell> 0$,
\[
\alpha(\ell) \eqdef\lim_{N \to\infty} \EE^{\pi^N} [ \alpha
^N(x,\xi^N) ]
= \EE[ 1 \wedge e^{Z_{\ell}}].
\]
This result is rigorously proved as Corollary~\ref{remconcentration}.
We then define the ``speed function''
\begin{equation} \label{espeedfunc}
h(\ell) = \ell\alpha(\ell).
\end{equation}
Note that the time step made in the proposal is $\delta=l\Delta t$
and that if this is accepted a fraction $\alpha(\ell)$ of the
time, then a naive argument invoking independence shows that the
effective time-step is reduced to $h(l)\Delta t.$
This is made rigorous in Theorem~\ref{thmmain}
which shows that the quantity $h(\ell)$ is the asymptotic speed function
of the limiting diffusion obtained by rescaling the Metropolis--Hastings
Markov chain $x^N=\{x^{k,N}\}_{k \geq0}$.

\begin{theorem}[(Main theorem)] \label{thmmain}
Let the reference measure $\pi_0$ and the function $\Psi(\cdot)$ satisfy
Assumption~\ref{ass1}. Consider the MALA algorithm
\eqref{eqnMALAalgorithm} with initial condition
$x^{0,N} \dist\pi^N.$ Let $z^N(t)$ be the piecewise linear,
continuous interpolant of the MALA algorithm
as defined in \eqref{eqnMCMCe}, with $\Delta t=N^{-1/3}.$
Then $z^N(t)$ converges weakly in $C([0,T],\h^s)$
to the diffusion process $z(t)$ given by
\begin{equation} \label{eqnspdemain}
\frac{dz}{dt} = - h(\ell) \bigl(z + \C\nabla\Psi(z) \bigr)+ \sqrt{2
h(\ell)} \,\frac{dW}{dt}
\end{equation}
with initial distribution $z(0) \dist\pi$.
\end{theorem}

We now explain the following two important implications of this result:
\begin{itemize}
\item
Since time has to be accelerated by a factor $(\Delta t)^{-1}
=N^{{1}/{3}}$
in order to observe a diffusion limit, it follows
that in stationarity the work required to
explore the invariant measure scales as $\OO(N^{{1}/{3}})$.
\item
The speed at which the invariant
measure is explored, again in stationarity, is maximized
by choosing $\ell$ so as to maximize $h(\ell)$;
this is achieved at an average acceptance probability $0.574$.
From a practical point of view, this shows that
one should ``tune'' the proposal variance of the MALA algorithm
so as to have a mean acceptance probability of $0.574$.
\end{itemize}

\begin{figure}

\includegraphics{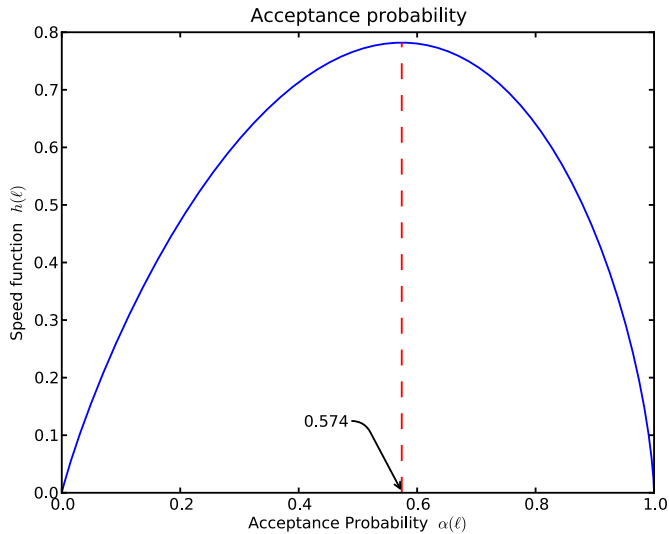}

\caption{Optimal acceptance probability${} = 0.574$.}\label{figoptaccept}
\end{figure}
\noindent The first implication follows from \eqref{eqnMCMCe} since this shows
that $\OO(N^{{1}/{3}})$ steps of the MALA Markov chain \eqref
{eqnMALAalgorithm} are
required for $z^N(t)$ to approximate $z(t)$ on a time interval $[0,T]$
long enough for $z(t)$ to have explored its invariant measure.
To understand the second implication, note that if $Z(t)$ solves
\eqref{eqnspdemain} with $h(\ell)\equiv1$, then, in law,
$z(t)=Z(h(\ell)t).$
This result suggests choosing the value of $\ell$ that maximizes
the speed function $h(\cdot)$ since $z(t)$ will then
explore the invariant measure as fast as possible.
For practitioners, who often tune algorithms according to
the acceptance probability, it is relevant to express the
maximization principle in terms of the asymptotic mean acceptance
probability $\alpha(\ell)$. Figure~\ref{figoptaccept} shows that the speed
function $h(\cdot)$ is maximized for an optimal acceptance probability
of $\alpha^{\star} =0.574$, to three-decimal places.
This is precisely the argument used in~\cite{RobeRose98}
for the case of product target measures, and it is remarkable
that the optimal acceptance probability identified in that context
is also optimal for the nonproduct measures studied in this paper.

\section{Proof of main theorem}
\label{secproof}

In Section~\ref{ssecps} we outline the proof strategy
and introduce the drift-martingale decomposition of our
discrete-time Markov chain which underlies it. Section~\ref{ssecGDA}
contains statement and proof of a general diffusion
approximation, Proposition~\ref{propdiffapprox}.
In Section~\ref{secproofmainthm} we use this
proposition to prove the main theorem of this
paper, pointing to Section~\ref{seckeyestimates}
for the key estimates required.

\subsection{Proof strategy}
\label{ssecps}
To communicate the main ideas, we give a heuristic of the proof
before proceeding to give full details in subsequent sections.
Let us first examine a simpler situation:
consider a scalar Lipschitz function $\mu\dvtx\RR\to\RR$ and two scalar
constants $\ell, c>0$. The usual theory
of diffusion approximation for Markov processes~\cite{EthiKurt86}
shows that the sequence $x^N=\{ x^{k,N} \}$ of Markov chains
\[
x^{k+1,N} - x^{k,N}= \mu(x^{k,N}) \ell N^{-1/3} + \sqrt{2 \ell
N^{-1/3}} c^{{1}/{2}} \xi^k,
\]
with i.i.d. $\xi^k \dist\Normal(0,1)$
converges weakly,
when interpolated using a time-acceleration factor
of $N^{1/3}$, to the scalar diffusion
$dz(t) = \ell\mu(z(t)) \,dt +  \sqrt{2\ell} \,dW(t)$ where
$W$ is a Brownian
motion with variance $\var(W(t))=ct$. Also, if $\gamma^k$ is
an $\iid$ sequence
of Bernoulli random variables with success rate $\alpha(\ell)$,
independent from the Markov chain $x^N$, one can prove
that the sequence $x^N=\{ x^{k,N} \}$ of Markov chains given by
\[
x^{k+1,N} - x^{k,N}= \gamma^k \bigl\{ \mu(x^{k,N}) \ell
N^{-1/3}+ \sqrt{2 \ell N^{-1/3}} c^{{1}/{2}} \xi^k \bigr\}
\]
converges weakly, when interpolated using a time-acceleration
factor $N^{1/3}$, to the diffusion
\[
dz(t) = h(\ell) \mu(z(t)) \,dt + \sqrt{2h(\ell)} \,dW(t),
\]
where the speed function is given by $h(\ell) = \ell\alpha(\ell)$.
This shows that the Bernoulli random variables $\{ \gamma^k \}
_{k \geq0}$ have slowed down the
original Markov chain by a factor $\alpha(\ell)$.
The proof of Theorem~\ref{thmmain} is an application of this idea in a
slightly more general setting. The following
complications arise:
\begin{itemize}
\item
Instead of working with scalar diffusions, the result holds for a
Hilbert space-valued diffusion.
The correlation structure between the different coordinates is not
present in the preceding
simple example and has to be taken into account.

\item
Instead of working with a single drift function $\mu$, a sequence of
approximations $d^N$ converging to $\mu$ has to be taken into account.

\item
The Bernoulli random variables $\gamma^{k,N}$ are not $\iid$ and
have an autocorrelation structure. On top of that, the Bernoulli random
variables $\gamma^{k,N}$ are not independent from the
Markov chain $x^{k,N}$. This is the main difficulty in the proof.

\item
It should be emphasized that the main theorem uses the fact that the
MALA Markov chain
is started at stationarity; this, in particular, implies that $x^{k,N}
\dist\pi^N$ for any $k \geq0$,
which is crucial to the proof of the invariance principle as
it allows us to control the correlation between $\gamma^{k,N}$
and $x^{k,N}$.

\end{itemize}

The acceptance probability of proposal \eqref{eqnproposal} is equal
to $\alpha^N(x,\xi^N) = 1 \wedge e^{Q^N(x,\xi^N)}$, and the quantity
$\alpha^N(x) = \EE_x[\alpha^N(x, \xi^N)]$, given by
\eqref{elocaccept}, represents the mean acceptance
probability when the Markov chain $x^N$ stands at $x$.
For our proof it is important to understand
how the acceptance probability $\alpha^N(x, \xi^N)$ depends on the current
position $x$ and on the source of randomness $\xi^N$.
Recall the quantity $Q^N$ defined in equation \eqref{eqnQN}: the main
observation
is that $Q^N(x,\xi^N)$ can be approximated by a Gaussian random variable
\begin{equation} \label{egaussapprox}
Q^N(x,\xi^N) \approx Z_{\ell},
\end{equation}
where $Z_{\ell} \dist\Normal(-\frac{\ell^3}{4}, \frac{\ell^3}{2})$.
These approximations
are made rigorous in Lem\-ma~\ref{lemgaussianapprox} and Lemma \ref
{lemconcentration}.
Therefore, the Bernoulli random variable $\gamma^{N}(x,\xi^N)$ with
success probability
$1 \wedge e^{Q^N(x,\xi^N)}$ can be approximated by a Bernoulli random
variable, independent of $x$,
with success probability equal to
\begin{equation} \label{elimitaccprob}
\alpha(\ell) = \EE[ 1 \wedge e^{Z_{\ell}} ].
\end{equation}
Thus, the limiting acceptance probability of the MALA algorithm
is as given in equation \eqref{elimitaccprob}.

Recall that $\Delta t=N^{-1/3}.$ With this notation
we introduce the drift function $d^N\dvtx\h^s \to\h^s$ given by
\begin{equation} \label{eqdrift}
d^N(x) = ( h(\ell) \Delta t)^{-1} \EE[ x^{1,N} - x^{0,N}
|x^{0,N} = x ]
\end{equation}
and the martingale difference array $\{\Gamma^{k,N}\dvtx k \geq0\}$
defined by $\Gamma^{k,N} = \Gamma^N(x^{k,N}, \break \xi^{k,N})$
with
\begin{equation} \label{eqmart}
\Gamma^{k,N} = ( 2 h(\ell) \Delta t )^{-1/2} \bigl(
x^{k+1,N} - x^{k,N} - h(\ell) \Delta t d^N(x^{k,N}) \bigr).
\end{equation}
The normalization constant $h(\ell)$ defined in equation \eqref{espeedfunc}
ensures that the drift function $d^N$ and the martingale difference
array $\{\Gamma^{k,N}\}$
are asymptotically independent from the parameter $\ell$.
The drift-martingale decomposition of the Markov chain $\{x^{k,N}\}_k$
then reads
\begin{equation} \label{eqndrift-mart}
x^{k+1,N} - x^{k,N} = h(\ell) \Delta t d^N(x^{k,N})
+ \sqrt{2 h(\ell) \Delta t} \Gamma^{k,N}.
\end{equation}
Lemma~\ref{lemdriftapprox} and Lemma~\ref{lemdiffus}
exploit the Gaussian behavior of $Q^N(x,\xi^N)$, described in equation
\eqref{egaussapprox},
in order to give quantitative versions of the following approximations:
\begin{equation} \label{edriftnoiseapprox}
d^N(x) \approx\mu(x)
\quad\mbox{and}\quad
\Gamma^{k,N} \approx\Normal(0,\C),
\end{equation}
where the function $\mu(\cdot)$ is defined by equation \eqref{eqnmu}.
From equation \eqref{eqndrift-mart} it follows that for large $N$
the evolution of the Markov chain resembles the Euler
discretization of the limiting diffusion \eqref{eqnspdemain}.
The next step consists of
proving an invariance principle for a rescaled version
of the martingale difference array $\{\Gamma^{k,N}\}$.
The continuous process $W^N \in C([0;T], \h^s)$
is defined as
\begin{equation} \label{eWN}
 W^N(t) = \sqrt{\Delta t} \sum_{j=0}^k \Gamma^{j,N}
+ \frac{t - k\Delta t}{\sqrt{\Delta t}} \Gamma^{k+1,N}
\qquad\mbox{for }
k \Delta t \leq t < (k+1) \Delta t.\hspace*{-35pt}
\end{equation}
The sequence of processes $\{W^N\}_{N \geq1}$ converges weakly as $N
\to\infty$ in
$C([0;T],\break  \h^s)$ to a Brownian motion $W$ in $\h^s$ with covariance
operator equal to $\C_s$.
Indeed, Proposition~\ref{lembweakconv} proves the stronger result
\[
(x^{0,N},W^N) \Longrightarrow(z^0,W),
\]
where $\Longrightarrow$ denotes weak convergence in $\h^s \times
C([0;T],\h^s)$, and $z^0 \dist\pi$ is independent
of the limiting Brownian motion $W$.
Using this invariance principle and the fact that
the noise process is additive [the diffusion coefficient
of the SPDE~\eqref{eqnspdemain} is constant],
the main theorem follows from a continuous mapping argument which we
now outline.
For any $W \in C([0,T]; \h^s)$ we define the It{\^o} map
\[
\Theta\dvtx \h^s \times C([0,T];\h^s) \rightarrow C([0,T];\h^s)
\]
which maps $(z^0, W)$ to the unique solution of the integral equation
\begin{equation}\label{eqnpcxN1jcm}
z(t) = z^0 - h(\ell) \int_0^t \mu(z) \,du + \sqrt{2h(\ell)} W(t)
\qquad
\forall t \in[0,T].
\end{equation}
Notice that $z=\Theta(z^0,W)$ solves the SPDE \eqref{eqnspdemain}.
The It{\^o} map $\Theta$ is continuous, essentially
because the noise in \eqref{eqnspdemain} is additive (does
not depend on the state $z$).
The piecewise constant interpolant $\bar{z}^N$ of $x^{N}$
is defined by
\begin{equation} \label{epicewisecstinterp}
\bar{z}^N(t) = x^k \qquad\mbox{for } k \Delta t \leq t < (k+1)
\Delta t.
\end{equation}
Using this definition it follows that
the continuous piecewise linear interpolant $z^N$, defined
in equation \eqref{eqnMCMCe}, satisfies
\begin{equation} \label{eqnprocessxN11}
z^N(t) = x^{0,N} - h(\ell) \int_0^t d^N(\bar{z}^N(u)) \,du + \sqrt
{2h(\ell)} W^N(t)\qquad
\forall t \in[0,T].\hspace*{-35pt}
\end{equation}
Using the closeness of $d^N(\cdot)$ and $\mu(\cdot)$,
and of $z^N$ and $\bar{z}^N$, we will see that there exists a process
$\widehat{W}^N \weak W$ as $N\rightarrow\infty$
such that
\[
z^N(t) = x^{0,N} - h(\ell) \int_0^t \mu( z^N(u) ) \,du + \sqrt
{2h(\ell)} \widehat W^N(t).
\]
Thus we may write $z^N=\Theta(x^{0,N}, \widehat W^N)$.
By continuity
of the It{\^o} map $\Theta$, it follows from the continuous mapping
theorem that
$ z^N = \Theta(x^{0,N}, \widehat W^N) \longweak\Theta(z^0, W) = z$
as $N$ goes to infinity. This weak convergence result is
the principal result of this article.

\subsection{General diffusion approximation}
\label{ssecGDA}
In this section we state and prove a proposition
containing a general diffusion approximation result.
Using this, we then prove our main theorem in Section~\ref{secproofmainthm}.
To this end, consider a general sequence of Markov chains $x^N=\{
x^{k,N}\}_{k \geq0}$
evolving at stationarity in the separable Hilbert space $\h^s$, and introduce
the drift-martingale decomposition
%
%
\begin{equation} \label{eqnmarkovchainseq}
x^{k+1,N}-x^{k,N} = h(\ell) d^N(x_k) \Delta t + \sqrt{2 h(\ell
) \Delta t} \Gamma^{k,N},
\end{equation}
where $h(\ell) > 0$ is a constant parameter, and $\Delta t$
is a time-step decreasing to $0$ as $N$ goes to infinity.
Here $d^N$ and $\Gamma^{k,N}$ are as defined above.
We introduce the rescaled process $W^N(t)$
as in \eqref{eWN}. The main diffusion approximation result
is the following.

\begin{prop}[(General diffusion approximation for Markov chains)] \label{propdiffapprox}
Consider a separable Hilbert space $( \h^s, \langle\cdot, \cdot
\rangle_s
)$
and a sequence of $\h^s$-valued Markov chains
$x^N = \{x^{k,N}\}_{k \geq0}$ with invariant
distribution $\pi^N$.
Suppose that the Markov chains start at stationarity $x^{0,N} \dist\pi^N$
and that the drift-martingale decomposition \eqref
{eqnmarkovchainseq} satisfies
the following assumptions:
\begin{longlist}[(1)]
\item[(1)]
\textit{Convergence of initial conditions}: $\pi^N$
converges in distribution to the probability measure $\pi$ where $\pi$
has a finite first moment, that is,
\mbox{$\EE^{\pi}[ \|x\|_s ] < \infty$}.
\item[(2)]
\textit{Invariance principle}:
the sequence $(x^{0,N}, W^N)$, defined by equation
\eqref{eWN}, converges weakly in
$\h^s \times C([0,T],\h^s)$ to $(z^0, W)$ where $z^0 \dist\pi$,
and $W$ is a Brownian motion in $\h^s$, independent from $z^0$,
with covariance operator~$\C_s$.
\item[(3)]
\textit{Convergence of the drift}:
There exists a globally Lipschitz function $\mu\dvtx\h^s \to\h^s$
that satisfies
\[
\lim_{N \to\infty} \EE^{\pi^N}[ \|d^N(x)-\mu(x)\|_s ] = 0.
\]
\end{longlist}
Then the sequence of rescaled interpolants
$z^N \in C([0,T],\h^s)$, defined by
equation \eqref{eqnMCMCe},
converges weakly in $C([0,T],\h^s)$
to $z \in C([0,T],\h^s)$ given by
%
\begin{eqnarray*}
\frac{dz}{dt} &=& h(\ell) \mu(z(t)) + \sqrt{2 h(\ell)}
\,\frac{dW}{dt},\\
z(0) &\dist&\pi.
\end{eqnarray*}
Here $W$ is a Brownian motion in $\h^s$ with covariance $\C_s$ and
initial condition $z^0 \dist\pi$ independent of $W$.
\end{prop}

\begin{pf}
Define $\bar{z}^N(t)$ as in \eqref{epicewisecstinterp}.
It then follows that
\begin{eqnarray} \label{einteq}
z^N(t) &=& x^{0,N} + h(\ell) \int_0^t d^N(\bar{z}^N(u)) \,du + \sqrt
{2h(\ell)} W^N(t)
\nonumber
\\[-8pt]
\\[-8pt]
\nonumber
&=& z^{0,N} + h(\ell) \int_0^t \mu(z^N(u)) \,du + \sqrt{2h(\ell)}
\widehat W^N(t),
\end{eqnarray}
where the process $W^N \in C([0,T], \h^s)$ is defined
by equation \eqref{eWN} and
\[
\widehat W^N(t) = W^N(t) + \sqrt{\frac{h(\ell)}{2}} \int_0^t [
d^N(\bar
{z}^N(u))- \mu(z^N(u)) ] \,du.
\]
Define the It{\^o} map
$\Theta\dvtx \h^s \times C([0,T];\h^s) \rightarrow C([0,T];\h^s)$
that maps $(z_0, W)$ to the unique solution $z \in C([0,T],\h^s)$
of the integral equation
\[
z(t) = z_0 + h(\ell) \int_0^t \mu(z(u)) \,du + \sqrt{2 h(\ell)}
W(t) \qquad\forall t \in[0,T].
\]
Equation \eqref{einteq} is thus equivalent to
$z^N = \Theta(x^{0,N}, \widehat W^N)$.
The proof of the diffusion approximation is accomplished through the
following steps:
\begin{itemize}
%
\item
\textit{The It{\^o} map $\Theta\dvtx\h^s \times C([0,T], \h^s) \to
C([0,T], \h
^s)$ is continuous}.
This is Lemma $3.7$ of~\cite{MatPillStu09}.
%
\item\textit{The pair $(x^{0,N},\widehat W^N)$ converges weakly to
$(z^0,W)$}.
In a separable Hilbert space, if the sequence $\{a_n\}_{n \in\NN}$
converges weakly to $a$, and
the sequence $\{b_n\}_{n \in\NN}$
converges in probability to $0$,
then the sequence $\{a_n+b_n\}_{n \in\NN}$
converges weakly to $a$.
It is assumed that $(x^{0,N}, W^N)$ converges weakly to $(z^0, W)$ in
$\h^s \times C([0,T], \h^s)$.
Consequently, to prove that $\widehat W^N$ converges weakly to $W$,
it suffices to prove that $\int_0^T \| d^N(\bar{z}^N(u)) - \mu(z^N(u))
\|_s \,du$ converges in probability to~$0$.
For any time $k \Delta t \leq u < (k+1) \Delta t$, the stationarity of
the chain shows that
\begin{eqnarray*}
\|d^N(\bar{z}^N(u)) - \mu(\bar{z}^N(u))\|_s &= &\|d^N(x^{k,N})
- \mu(x^{k,N})\|_s\\
& \dist&\|d^N(x^{0,N}) - \mu(x^{0,N})\|_s, \\
\|\mu(\bar{z}^N(u)) - \mu(z^N(u))\|_s &\leq&\|\mu\|_{\mathrm{Lip}}
\cdot
\|x^{k+1,N} - x^{k,N}\|_s\\
& \dist&\|\mu\|_{\mathrm{Lip}} \cdot\|x^{1,N}
- x^{0,N}\|_s,
\end{eqnarray*}
where in the last step we have used the fact that $\|\bar{z}^N(u) -
z^N(u)\|_s \leq\|x^{k+1,N} - x^{k,N}\|_s$. Consequently,
\begin{eqnarray*}
&&\EE^{\pi^N}\biggl[ \int_0^T \| d^N(\bar{z}^N(u)) - \mu(z^N(u)) \|_s \,du
\biggr] \\
&&\qquad\leq
T \cdot\EE^{\pi^N} [ \|d^N(x^{0,N}) - \mu(x^{0,N})\|_s ]\\
&&\qquad\quad{}+ T \cdot\|\mu\|_{\mathrm{Lip}} \cdot\EE^{\pi^N} [ \|x^{1,N} -
x^{0,N}\|_s ].
\end{eqnarray*}
The first term goes to zero since it is assumed that $\lim_N \EE^{\pi
^N}[ \| d^N(x)-\mu(x)\|_s] = 0$.
Since $\tr_{\h^s}(\C_s) < \infty$,
the second term is of order $\OO(\sqrt{\Delta t})$ and thus also
converges to $0$.
Therefore $\widehat{W}^N$ converges weakly to $W$, hence the conclusion.

%
%
\item\textit{Continuous mapping argument}.
We have proved that $(x^{0,N}, \widehat W^N)$ converges weakly
in $\h^s \times C([0,T],\h^s)$ to $(z^0,W)$, and the It{\^o} map
$\Theta\dvtx \h^s \times C([0,T], \h^s) \to C([0,T],\h^s)$ is a
continuous function. The continuous mapping theorem thus shows that
$z^N = \Theta(x^{0,N}, \widehat W^N)$
converges weakly to $z=\Theta(z^0,W)$, finishing the proof of
Proposition~\ref{propdiffapprox}.\quad\qed
\end{itemize}
\noqed\end{pf}

\subsection{Proof of main theorem}
\label{secproofmainthm}

We now prove Theorem~\ref{thmmain}.
The proof consists of checking that
the conditions needed for Proposition~\ref{propdiffapprox} to apply
are satisfied by the sequence of MALA Markov
chains \eqref{eqnMALAalgorithm}. The key estimates are proved later in
Section~\ref{seckeyestimates}.
\begin{longlist}[(1)]
\item[(1)]
By Lemma~\ref{lemmchmeas} the sequence of probability measures $\pi^N$
converges weakly in $\h^s$ to $\pi$.
\item[(2)]
Proposition~\ref{lembweakconv} proves that
$(x^{0,N}, W^N)$ converges weakly in $\h\times   C([0,T], \h^s)$ to
$(z^0, W)$,
where $W$ is a Brownian motion with covariance $\C_s$ independent from
$z^0 \dist\pi.$

\item[(3)]
Lemma~\ref{lemdriftapprox} states that $d^N(x)$, defined by equation
\eqref{eqdrift}, satisfies
$\lim_N \EE^{\pi^N} [ \|d^N(x)-\mu(x)\|_s^2 ] = 0$, and
Proposition~\ref{lemregularity}
shows that $\mu\dvtx\h^s \to\h^s$ is a Lipschitz function.
\end{longlist}
The three assumptions needed for Lemma~\ref{propdiffapprox} to
apply are satisfied, which concludes the proof of Theorem~\ref{thmmain}.

\section{Key estimates}
\label{seckeyestimates}

Section~\ref{ssectl} contains some technical
lemmas of use \mbox{throughout}. In Section~\ref{ssecga}
we study the large $N$ Gaussian approximation of the acceptance
probability, simultaneously establishing asymptotic independence
of the current state of the Markov chain. This
approximation is then used in Sections~\ref{ssecda}
and~\ref{ssecna} to give quantitative versions of the
heuristics \eqref{edriftnoiseapprox}.
The section concludes with Section~\ref{secbweakconv}
in which we prove an invariance principle for $W^N$
given by \eqref{eWN}.

\subsection{Technical lemmas}
\label{ssectl}

The first lemma shows that, for $\pi_0$-almost every function $x \in
\h^s$,
the approximation $\mu^N(x) \approx\mu(x)$ holds as $N$ goes to infinity.

\begin{lemma}[($\mu^N$ converges $\pi_0$-almost
surely to $\mu$)] \label{lemmuNmu}
Let Assumption~\ref{ass1} hold. The sequences of functions $\mu
^N\dvtx\h^s
\to\h^s$ satisfies
\[
\pi_0 \Bigl( \Bigl\{ x \in\h^s\dvtx\lim_{N \to\infty} \| \mu^N(x) -
\mu(x) \|_s = 0 \Bigr\} \Bigr) = 1.
\]
\end{lemma}

\begin{pf}
It is enough to verify that for $x \in\h^s$, we have
\begin{eqnarray}
\lim_{N \to\infty} \|P^N x - x \|_s &= &0, \label{ealmostsure1} \\
\lim_{N \to\infty} \|\C P^N \nabla\Psi( P^N x) - \C\nabla\Psi
(x) \|
_s& =& 0 \label{ealmostsure2}.
\end{eqnarray}
\begin{itemize}
\item
Let us prove equation \eqref{ealmostsure1}.
For $x \in\h^s$, we have $\sum_{j \geq1} j^{2s} x_j^2 < \infty$ so that
\begin{equation} \label{eCDT}
\lim_{N \to\infty} \|P^N x - x \|_s^2
= \lim_{N \to\infty} \sum_{j=N+1}^{\infty} j^{2s} x_j^2 = 0.
\end{equation}
\item Let us prove \eqref{ealmostsure2}. The triangle inequality
shows that
\begin{eqnarray*}
&&\|\C P^N \nabla\Psi( P^N x) - \C\nabla\Psi(x) \|_s\\
&&\qquad\leq\|\C P^N \nabla\Psi( P^N x) - \C P^N \nabla\Psi(x)\|_s
+ \|\C P^N \nabla\Psi(x) - \C\nabla\Psi(x)\|_s.
\end{eqnarray*}
The same proof as Lemma~\ref{lemregularity} reveals that $\C P^N
\nabla\Psi\dvtx\h^s \to\h^s$
is globally Lipschitz, with a Lipschitz constant that can be chosen\vadjust{\goodbreak}
independently of~$N$.
Consequently, equation \eqref{eCDT} shows that
\[
\|\C P^N \nabla\Psi( P^N x) - \C P^N \nabla\Psi(x)\|_s \lesssim\|P^N
x - x \|_s \to0.
\]
Also, $z = \nabla\Psi(x) \in\h^{-s}$ so that $\| \nabla\Psi(x)\|
_{-s}^2 = \sum_{j \geq1} j^{-2s} z_j^2 < \infty$.
The eigenvalues of $\C$ satisfy $\lambda_j^2 \asymp j^{-2\kappa}$ with
$s < \kappa- \frac{1}{2}$.
Consequently,
\begin{eqnarray*}
&&\|\C P^N \nabla\Psi(x) - \C\nabla\Psi(x)\|_s^2\\
&&\qquad= \sum_{j=N+1}^{\infty} j^{2s} (\lambda_j^2 z_j)^2
\lesssim\sum_{j=N+1}^{\infty} j^{2s-4\kappa} z_j^2 \\
&&\qquad= \sum_{j=N+1}^{\infty} j^{4(s-\kappa)} j^{-2s} z_j^2
\leq\frac{1}{(N+1)^{4(\kappa-s)}} \| \nabla\Psi(x)\|_{-s}^2 \to0.
\end{eqnarray*}
\end{itemize}
\upqed\end{pf}

The next lemma shows that the size of the jump $y-x$ is of order $\sqrt
{\Delta t}$.

\begin{lemma} \label{lemsizeprop}
Consider $y$ given by \eqref{eqnproposal}.
Under Assumption~\ref{ass1}, for any \mbox{$p\geq1$},
we have
\[
\EE_x^{\pi^N} [ \|y-x\|^p_s ] \lesssim(\Delta t)^{{p}/{2}}
\cdot(1 + \|x\|^p_s).
\]
\end{lemma}

\begin{pf}
Under Assumption~\ref{ass1} the function $\mu^N$ is globally Lipschitz
on $\h^s$, with Lipschitz constant that can be chosen independently of $N$.
Thus
\[
\|y-x\|_s \lesssim\Delta t (1+\|x\|_s) + \sqrt{\Delta t} \| \C
^{{1}/{2}}\xi^N \|_s .
\]
We have $\EE^{\pi^0} [ \|\C^{{1}/{2}} \xi^N\|^p_s ] \leq\EE
^{\pi^0} [ \|\zeta\|^p_s ] < \infty$,
where $\zeta\dist\Normal(0,\C)$. From Fernique's theorem~\cite{DaprZaby92},
it follows that $\EE^{\pi^0} [ \|\zeta\|^p_s ] < \infty$.
Consequently,\break $\EE^{\pi^0} [ \|\C^{{1}/{2}} \xi^N\|^p_s ]$
is uniformly bounded as a function of $N$,
proving the lemma.
\end{pf}

The normalizing constants $M_{\Psi^N}$ are uniformly bounded,
and we use this fact to obtain uniform bounds on moments of functionals
in~$\h$ under $\pi^N$. Moreover,
we prove that the sequence of probability measures $\pi^N$ on $\h^s$
converges weakly in~$\h^s$ to $\pi$.
%
\begin{lemma}[(Finite dimensional
approximation $\pi^N$ of $\pi$)] \label{lemmchmeas}
Under Assumption~\ref{ass1} the normalization
constants $M_{\Psi^N}$ are uniformly bounded so that for any measurable
functional $f\dvtx\h\mapsto\mathbb{R}$,
we have
\[
\EE^{\pi^N} [ |f(x)| ] \lesssim\EE^{\pi_0}[ |f(x)| ].
\]
Moreover, the sequence of probability measure $\pi^N$ satisfies
\[
\pi^N \longweak\pi,
\]
where $\longweak$ denotes weak convergence in $\h^s$.
\end{lemma}

\begin{pf}
The first part is contained in Lemma $3.5$ of~\cite{MatPillStu09}.
Let us prove that $\pi^N \longweak\pi$.
We need to show that for any bounded continuous function $g\dvtx\h^s
\to\RR
$ we have
$\lim_{N \to\infty} \EE^{\pi^N}[ g(x) ] = \EE^{\pi}[ g(x) ]$ where
\begin{eqnarray*}
\EE^{\pi^N}[ g(x) ]
&=& \EE^{\pi_0^N}\bigl[g(x) M_{\Psi^N} e^{-\Psi^N(x)} \bigr] \\
&=& \EE^{\pi_0}\bigl[ g( P^N x) M_{\Psi^N} e^{-\Psi( P^N x)} \bigr].
\end{eqnarray*}
Since $g$ is bounded, $\Psi$ is lower bounded, and since
the normalization constants are uniformly bounded,
the dominated convergence theorem shows that it suffices to show
that $g( P^N x) M_{\Psi^N} e^{-\Psi( P^N x)}$ converges $\pi_0$-almost
surely to $g(x) M_{\Psi} e^{-\Psi(x)}$. For this in turn
it suffices to show that $\Psi( P^N x)$
converges $\pi_0$-almost surely to $\Psi(x)$, as this also
proves almost sure convergence of the normalization constants.
By \eqref{eqnpsi2} we have
\[
|\Psi( P^N x)-\Psi(x)| \lesssim(1+\|x\|_s+\|P^N x\|_s)
\|P^N x - x\|_s.
\]
But $\lim_{N \to\infty} \|P^N x - x\|_s \to0$
for any $x \in\h^s$, by dominated convergence, and the
result follows.
\end{pf}
Fernique's theorem~\cite{DaprZaby92} states that for any exponent $p
\geq0$,
we have $\EE^{\pi^0}  [ \|x\|_s^p ] < \infty$. It thus follows from
Lemma~\ref{lemmchmeas} that for any $p \geq0$,
\[
\sup_N \{ \EE^{\pi^N}[ \|x\|_s^p ] \dvtx N \in\mathbb{N}
\} < \infty.
\]
This estimate is repeatedly used in the sequel.

\subsection{Gaussian approximation of $Q^N$}
\label{ssecga}

Recall the quantity $Q^N$ defined in equation \eqref{eqnQN}. This
section proves that $Q^N$
has a Gaussian behavior in the sense that
\begin{equation} \label{egaussapperror}
Q^N(x, \xi^N) = Z^N(x, \xi^N) + i^N(x,\xi^N) + \err^N(x,\xi^N),
\end{equation}
where the quantities $Z^N$ and $i^N$ are equal to
\begin{eqnarray}
Z^N(x, \xi^N)&=& -\frac{\ell^3}{4} - \frac{\ell^{{3}/{2}}}{\sqrt{2}}
N^{-{1}/{2}} \sum_{j=1}^N \lambda_j^{-1} \xi_j x_j, \label
{eZn}\\
i^N(x,\xi^N) &=& \frac{1}{2} (\ell\Delta t)^2 \bigl( \|x\|^2_{\C^N} - \|
(\C^N)^{{1}/{2}}\xi^N\|^2_{\C^N} \bigr) \label{eiN}
\end{eqnarray}
with $i^N$ and $e^N$ small.
Thus the principal contributions to $Q^N$
comes from the random variable $Z^N(x,\xi^N)$.
Notice that, for each fixed $x \in\h^s$, the
random variable $Z^N(x,\xi^N)$ is Gaussian.
Furthermore, the Karhunen--Lo\`eve expansion of
$\pi_0$ shows that for $\pi_0$-almost every choice of
function $x \in\h$ the sequence
$\{ Z^N(x, \xi^N) \}_{N \geq1}$ converges
in law to the distribution of
$Z_{\ell} \dist\Normal(-\frac{\ell^3}{4}, \frac{\ell^3}{2})$.
The next lemma rigorously bounds\vadjust{\goodbreak} the error terms
$\err^N(x,\xi^N)$ and $i^N(x,\xi^N)$: we show that $i^N$ is an error
term of order $\OO(N^{-1/6})$
and $\err^N(x,\xi)$ is an error term of order $\OO(N^{-1/3})$.
In Lemma~\ref{lemconcentration} we then quantify the convergence
of $Z^N(x,\xi^N)$ to $Z_{\ell}.$

\begin{lemma}[(Gaussian approximation)] \label{lemgaussianapprox}
Let $p \geq1$ be an integer.
Under Assumption~\ref{ass1}, the error terms
$i^N$ and $\err^N$ in the Gaussian approximation \eqref
{egaussapperror} satisfy
%
%
\begin{eqnarray} \label{eqnerrorcontrol}
(\EE^{\pi^N} [|i^N(x,\xi^N)|^p ])^{{1}/{p}} &=& \OO
(N^{-{1}/{6}})
\quad\mbox{and}
\nonumber
\\[-8pt]
\\[-8pt]
\nonumber
(\EE^{\pi^N} [ |\err^N(x,\xi^N)|^p ] )^{{1}/{p}} &=&
\OO(N^{-{1}/{3}}).
\end{eqnarray}
\end{lemma}

\begin{pf}
For notational clarity, without loss of generality, we suppose $p=2q$.
The quantity $Q^N$ is defined in equation \eqref{eqnQN}, and
expanding terms leads to
\[
Q^N(x,\xi^N) = I_1 + I_2 + I_3,
\]
where the quantities $I_1$, $I_2$ and $I_3$ are given by
\begin{eqnarray*}
I_1 &=& -\frac{1}{2} ( \|y\|_{\C^N}^2 - \|x\|_{\C^N}^2 )\\
&&{}- \frac{1}{4 \ell\Delta t} \bigl( \|x-y(1-\ell\Delta t)\|_{\C^N}^2 - \|
y-x(1-\ell\Delta t)\|_{\C^N}^2 \bigr), \\
I_2 &= &-\bigl(\Psi^N(y) - \Psi^N(x) \bigr)
- \frac{1}{2} \bigl( \langle x-y(1-\ell\Delta t), \C^N \nabla\Psi ^N(y)
\rangle_{\C^N}\\
&&\hspace*{114pt}{}- \langle y-x(1-\ell\Delta t), \C^N \nabla\Psi^N(x) \rangle_{\C
^N}\bigr), \\
I_3 &=& -\frac{\ell\Delta t}{4}\{ \|\C^N \nabla\Psi^N(y)\|_{\C
^N}^2 - \|\C^N \nabla\Psi^N(x)\|_{\C^N}^2\}.
\end{eqnarray*}
The term $I_1$ arises purely from the Gaussian part of the
target measure $\pi^N$ and from the Gaussian part of
the proposal. The two other terms $I_2$
and $I_3$ come from the change of probability involving the functional
$\Psi^N$.
We start by simplifying the expression for $I_1$, and then
return to estimate the terms $I_2$ and $I_3$:
\begin{eqnarray*}
I_1
&=& -\frac{1}{2} ( \|y\|_{\C^N}^2 - \|x\|_{\C^N}^2 )
\\
&&{}- \frac{1}{4 \ell\Delta t} \bigl( \|(x-y) + \ell\Delta t y\|_{\C
^N}^2 - \|(y-x) + \ell\Delta t x\|_{\C^N}^2 \bigr) \\
&=& -\frac{1}{2} ( \|y\|_{\C^N}^2 - \|x\|_{\C^N}^2 )
\\
&&{}- \frac{1}{4 \ell\Delta t} \bigl( 2 \ell\Delta t [\|x\|_{\C^N}^2 - \|
y\|_{\C^N}^2] + (\ell\Delta t)^2 [\|y\|_{\C^N}^2 - \|x\|_{\C^N}^2]
\bigr) \\
&=& -\frac{\ell\Delta t}{4} ( \|y\|_{\C^N}^2 - \|x\|_{\C^N}^2).
\end{eqnarray*}
The term $I_1$ is $\OO(1)$ and constitutes the main contribution to
$Q^N$. Before analyzing~$I_1$
in more detail, we show that $I_2$ and $I_3$ are $\OO(N^{-{1}/{3}})$.
\begin{equation} \label{eI23small}
\qquad(\EE^{\pi^N}[I_2^{2q}])^{{1}/{(2q)}} = \OO(N^{-1/3})
\quad\mbox{and}\quad
(\EE^{\pi^N}[I_3^{2q}])^{{1}/{(2q)}} = \OO(N^{-1/3}).
\end{equation}
\begin{itemize}
\item
We expand $I_2$ and use the bound on the remainder of the Taylor expansion
of $\Psi$ described in equation \eqref{eqn2nd-Taylor},
\begin{eqnarray*}
I_2
&=& - \{ \Psi^N(y) - [\Psi^N(x) + \langle\nabla\Psi^N(x), y-x
\rangle]
\}\\
&&{} + \frac{1}{2}\langle y-x, \nabla\Psi^N(y)-\nabla\Psi^N(x)
\rangle \\
& &{}+ \frac{\ell\Delta t}{2} \{ \langle x, \nabla\Psi^{N}(x) \rangle -
\langle y, \nabla\Psi^{N}(y) \rangle\} \\
&=& A_1 + A_2 + A_3.
\end{eqnarray*}
Equation \eqref{eqn2nd-Taylor} and Lemma~\ref{lemsizeprop} show that
\[
\EE^{\pi^N}[ A_1^{2q} ] \lesssim\EE^{\pi^N}[ \|y-x\|_s^{4q} ]
\lesssim
(\Delta t)^{2q} \EE^{\pi^N}[1+\|x\|_s^{4q}]
\lesssim(\Delta t)^{2q} = ( N^{-1/3} )^{2q},
\]
where we have used the fact that $\EE^{\pi^N}[\|x\|_s^{4q}] \lesssim
\EE
^{\pi_0}[\|x\|_s^{4q}] < \infty$.
Assumption~\ref{ass1} states that $\partial^2 \Psi$ is uniformly
bounded in $\mathcal{L}(\h^s,\h^{-s})$ so that
\begin{eqnarray} \label{epsis}
\|\nabla\Psi(y) - \nabla\Psi(y)\|_{-s}
&=& \biggl\| \int_{0}^1 \partial^2 \Psi\bigl(x+t(y-x) \bigr) \cdot(y-x)
\,dt \biggr\|_{-s} \nonumber\\
&\leq&\int_{0}^1 \bigl\| \partial^2 \Psi\bigl(x+t(y-x) \bigr) \cdot(y-x) \bigr\|
_{-s} \,dt\\
&\leq& M_4 \int_{0}^1 \| y-x \|_s \,dt.\nonumber
\end{eqnarray}
This proves that $\|\nabla\Psi^N(y)-\nabla\Psi^N(x)\|_{-s} \lesssim
\|
y-x\|_s$.
Consequently, Lem\-ma~\ref{lemsizeprop} shows that
\begin{eqnarray*}
\EE^{\pi^N}[A_2^{2q}]
&\lesssim&\EE^{\pi^N}[ \|y-x\|_s^{2q} \cdot\|\nabla\Psi
^N(y)-\nabla\Psi^N(x)\|_{-s}^{2q} ]\\
&\lesssim&\EE^{\pi^N}[ \|y-x\|_s^{4q} ] \\
& \lesssim&(\Delta t)^{2q} \EE^{\pi^N} [1+\|x\|_s^{4q} ]\\
&\lesssim&(\Delta t)^2 = (N^{-1/3} )^{2q}.
\end{eqnarray*}
Under Assumption~\ref{ass1}, for any $z\in\h^s$ we have $\|\nabla
\Psi
^N(z)\|_{-s} \lesssim1 + \|z\|_s$.
Therefore $\EE^{\pi^N}[A_3^{2q}] \lesssim(\Delta t)^{2q}$. Putting
these estimates together,
\[
( \EE^{\pi^N}[I_2^{2q}] )^{{1}/{(2q)}} \lesssim( \EE
^{\pi^N}[A_1^{2q} + A_2^{2q} + A_3^{2q}] )^{{1}/{(2q)}} = \OO
(N^{-1/3}).
\]

\item
Lemma~\ref{lemregularity} states $\C^N \nabla\Psi^N\dvtx\h^s \to
\h^s$ is
globally Lipschitz,
with a Lipschitz constant that can be chosen uniformly in $N$. Therefore,
\begin{equation} \label{eCPsibound}
\| \C^N \nabla\Psi^N(z) \|_s \lesssim1 + \| z \|_s.
\end{equation}
Since $\|\C^N \nabla\Psi^N(z)\|_{\C^N}^2 = \langle\nabla\Psi
^N(z), \C^N \nabla\Psi^N(z) \rangle$, bound \eqref{eqnpsi2} gives
\begin{eqnarray*}
\EE^{\pi^N}[ I_3^{2q} ]
&\lesssim&\Delta t^{2q} \EE[ \langle\nabla\Psi^N(x), \C^N \nabla
\Psi^N(x) \rangle^q + \langle\nabla\Psi^N(y), \C^N \nabla\Psi
^N(y) \rangle^q ] \\
&\lesssim&\Delta t^{2q} \EE^{\pi^N}[ (1+\|x\|_s)^{2q} + (1+\|y\|
_s)^{2q} ] \\
&\lesssim&\Delta t^{2q} \EE^{\pi^N}[ 1+\|x\|^{2q}_s + \|y\|
^{2q}_s ] \lesssim\Delta t^{2q}
= (N^{-1/3})^{2q},
\end{eqnarray*}
which concludes the proof of equation \eqref{eI23small}.
\end{itemize}

We now simplify further the expression for $I_1$ and demonstrate that
it has a Gaussian behavior.
We use the definition of the proposal $y$, given in equation~\eqref
{eqnproposal}, to expand $I_1$.
For $x \in X^N$ we have $P^N x = x$. Therefore, for $x \in X^N$,
\begin{eqnarray*}
I_1
&=& -\frac{\ell\Delta t}{4} \bigl( \bigl\|(1-\ell\Delta t)x - \ell\Delta t
\C^N \nabla\Psi^N(x)
+ \sqrt{2 \ell\Delta t} (\C^N)^{{1}/{2}}\xi^N\bigr\|_{\C^N}^2 - \|
x\|_{\C^N}^2\bigr) \\
&=& Z^N(x,\xi^N) + i^N(x,\xi^N) + B_1 + B_2 + B_3 + B_4,
\end{eqnarray*}
with $Z^N(x,\xi^N)$ and $i^N(x,\xi^N)$ given by equation \eqref{eZn}
and \eqref{eiN} and
\begin{eqnarray*}
%
B_1 &= &\frac{\ell^3}{4} \biggl(1 - \frac{\|x\|_{\C^N}^2}{N}\biggr),\\
B_2 &=& -\frac{\ell^3}{4} N^{-1} \{ \|\C^N \nabla\Psi^N(x)\|^2_{\C
^N} + 2\langle x,\nabla\Psi^N(x) \rangle \} ,\\
B_3& =& \frac{\ell^{5/2}}{\sqrt{2}} N^{-5/6}\langle x+\C^N
\nabla \Psi^N(x), (\C^N)^{1/2} \xi^N \rangle_{\C^N},\\
B_4 &=& \frac{\ell^2}{2} N^{-2/3} \langle x, \nabla\Psi^N(x)
\rangle.
\end{eqnarray*}
The quantity $Z^N$ is the leading term. For each fixed value of $x \in
\h^s$,
the term $Z^N(x,\xi^N)$ is Gaussian. Below, we prove that quantity
$i^N$ is $\OO(N^{-1/6})$.
We now establish that each $B_j$ is $\OO(N^{-1/3})$,
\begin{equation} \label{eboundBj}
(\EE^{\pi^N} [ B_j^{2q} ] )^{{1}/{(2q)}} = \OO
(N^{-1/3})\qquad
j=1, \ldots, 4.
\end{equation}
\begin{itemize}
\item
Lemma~\ref{lemmchmeas} shows that
$\EE^{\pi^N}[(1 - \frac{\|x\|_{\C^N}^2}{N})^{2q}] \lesssim\EE
^{\pi_0}[(1 - \frac{\|x\|_{\C^N}^2}{N})^{2q}]$.
Under $\pi_0$,
\[
\frac{\|x\|_{\C^N}^2}{N} \dist\frac{\rho_1^2 + \cdots+
\rho^2_N}{N},
\]
where $\rho_1, \ldots, \rho_N$ are i.i.d. $\Normal(0,1)$ Gaussian random
variables.
Consequently, $\EE^{\pi^N}[B_1^{2q}]^{{1}/{(2q)}} = \OO(N^{-{1}/{2}})$.

\item
The term $\|\C^N \nabla\Psi^N(x)\|^{2q}_{\C^N}$ has already been
bounded while proving\break $\EE^{\pi^N} [I_3^{2q}] \lesssim( N^{-1/3} )^{2q}$.
Equation \eqref{eqnpsi2} gives the bound
$\|\nabla\Psi^N(x)\|_{-s} \lesssim1+\|x\|_s$ and shows that
$\EE^{\pi^N}[ \langle x,\nabla\Psi^N(x) \rangle^{2q} ]$ is uniformly
bounded as a function of $N$. Consequently,
\[
\EE^{\pi^N} [B_2^{2q} ]^{{1}/{(2q)}} = \OO(N^{-1}).
\]

\item
We have $\langle\C^N \nabla\Psi^N(x), (\C^N)^{1/2} \xi^N
\rangle_{\C
^N} =
\langle\nabla\Psi^N(x), (\C^N)^{1/2} \xi^N \rangle$ so that
\begin{eqnarray*}
\EE^{\pi^N}[\langle\C^N \nabla\Psi^N(x), (\C^N)^{1/2} \xi
^N \rangle_{\C^N}^{2q}]
&\lesssim&\EE^{\pi^N}[\|\nabla\Psi^N(x)\|_{-s}^{2q} \cdot\|(\C
^N)^{1/2} \xi^N\|_s^{2q}]\\
& \lesssim&1.
\end{eqnarray*}
By Lemma~\ref{lemmchmeas}, one can suppose $x \dist\pi_0$,
\[
\langle x, (\C^N)^{1/2} \xi^N \rangle_{\C^N} \dist\sum
_{j=1}^N \rho_j
\xi_j,
\]
where $\rho_1, \ldots, \rho_N$ are i.i.d. $\Normal(0,1)$ Gaussian random
variables.
Consequently $(\EE^{\pi^N}[\langle x, (\C^N)^{1/2} \xi^N
\rangle_{\C
^N}^{2q} ])^{{1}/{(2q)}} = \OO(N^{1/2})$, which proves that
\[
( \EE^{\pi^N}[ B_3^{2q} ] )^{{1}/{(2q)}} = \OO
(N^{-5/6 + 1/2}) = \OO(N^{-1/3}).
\]

\item
The bound $\|\nabla\Psi^N(x)\|_{-s} \lesssim1+\|x\|_s$ ensures that
$(\EE^{\pi^N} [ B_4^{2q} ])^{{1}/{(2q)}} =  \OO
(N^{-2/3})$.
\end{itemize}

Define the quantity $\err^N(x,\xi^N) = I_2 + I_3 + B_1+B_2+B_3+B_4$ so
that $Q^N$ can also be expressed as
\[
Q^N(x,\xi^N) = Z^N(x,\xi^N) + i^N(x,\xi^N) + \err^N(x,\xi^N).
\]
Equations \eqref{eI23small} and \eqref{eboundBj} show that $\err
^N$ satisfies
\[
( \EE^{\pi^N} [ \err^N(x,\xi^N)^{2q} ] )^{{1}/{(2q)}}
= \OO(N^{-1/3}).
\]
We now prove that $i^N$ is $\OO(N^{-1/6})$. By Lemma~\ref{lemmchmeas},
$\EE^{\pi^N}[i^N(x,\xi^N)^{2q}] \lesssim\EE^{\pi_0}[i^N(x,\xi^N)^{2q}]$.
If $x \dist\pi_0$, we have
\begin{eqnarray*}
i^N(x,\xi^N)
&= &\frac{\ell^2}{2} N^{-2/3}\{ \|x\|^2_{\C^N} - \|(\C
^N)^{1/2} \xi^N\|^2_{\C^N}\} \\
&=& \frac{\ell^2}{2} N^{-2/3} \sum_{j=1}^N (\rho_j^2 - \xi_j^2),
\end{eqnarray*}
where $\rho_1, \ldots, \rho_N$ are i.i.d. $\Normal(0,1)$ Gaussian random
variables.
Since\break $\EE[ \{\sum_{j=1}^N (\rho_j^2 - \xi_j^2)\}^{2q}
] \lesssim N^q$, it follows that
\begin{equation} \label{einbound}
( \EE^{\pi^N} [ i^N(x,\xi^N)^{2q} ] )^{{1}/{(2q)}}
= \OO(N^{-2/3 + 1/2}) = \OO(N^{-1/6}),
\end{equation}
which ends the proof of Lemma~\ref{lemgaussianapprox}.\vadjust{\goodbreak}
\end{pf}

The next lemma quantifies the fact that $Z^N(x,\xi^N)$ is
asymptotically independent from the current position $x$.

\begin{lemma}[(Asymptotic independence)] \label{lemconcentration}
Let $p \geq1$ be a positive integer and $f\dvtx\RR\to\RR$ be a
$1$-Lipschitz function.
Consider error terms $\err^N_{\star}(x,\xi)$ satisfying
\[
\lim_{N \to\infty} \EE^{\pi^N}[\err^N_{\star}(x,\xi^N)^p] = 0.
\]
Define the functions $\bar{f}^N\dvtx\RR\to\RR$ and the constant
$\bar{f}
\in\RR$ by
\[
\bar{f}^N(x) = \EE_x\bigl[ f\bigl(Z^N(x,\xi^N) + \err^N_{\star}(x,\xi^N)
\bigr) \bigr]\quad
\mbox{and}\quad
\bar{f} = \EE[f(Z_{\ell})].
\]
Then the function $f^N$ is highly concentrated around its mean in the
sense that
\[
\lim_{N \to\infty} \EE^{\pi^N}[ |\bar{f}^N(x) - \bar{f}|^p
] = 0.
\]
\end{lemma}

\begin{pf}
Let $f$ be a $1$-Lipschitz function. Define the function $F\dvtx\RR
\times
[0;\infty) \to\RR$ by
\[
F(\mu,\sigma) = \EE[ f(\rho_{\mu,\sigma})]\qquad
\mbox{where }
\rho_{\mu,\sigma} \dist\Normal(\mu, \sigma^2).
\]
The function $F$ satisfies
\begin{equation} \label{eFboun}
| F(\mu_1,\sigma_1) - F(\mu_2,\sigma_2)|
\lesssim
|\mu_2-\mu_1| + |\sigma_2-\sigma_1|,
\end{equation}
for any choice $\mu_1, \mu_2 \in\RR$ and $\sigma_1, \sigma_2 \geq0$.
Indeed,
\begin{eqnarray*}
| F(\mu_1,\sigma_1) - F(\mu_2,\sigma_2)|
&= &| \EE[ f(\mu_1 + \sigma_1 \rho_{0,1}) - f(\mu_2 +
\sigma_2 \rho_{0,1}) ] |\\
&\leq&\EE[ |\mu_2-\mu_1| + |\sigma_2-\sigma_1| \cdot|\rho
_{0,1}| ]\\
&\lesssim&|\mu_2-\mu_1| + |\sigma_2-\sigma_1|.
\end{eqnarray*}
We have $\EE_x[Z^N(x,\xi^N)] = \EE[Z_{\ell}] = -\frac{\ell^3}{4}$ while
the variances are given by
\[
\var[ Z^N(x,\xi^N) ] = \frac{\ell^3}{2} \frac{\|x\|^2_{\C^N}}{N}
\quad\mbox{and}\quad
\var[ Z_{\ell} ] = \frac{\ell^3}{2}.
\]
Therefore, using Lemma~\ref{lemmchmeas},
\begin{eqnarray*}
&&\EE^{\pi^N}[ | \bar{f}^N(x)-\bar{f} |^p ]\\
&&\qquad= \EE^{\pi^N}\bigl[ \bigl| \EE_x\bigl[ f\bigl(Z^N(x,\xi^N) + \err^N_{\star
}(x,\xi^N) \bigr) -f(Z_{\ell}) \bigr] \bigr|^p \bigr]\\
&&\qquad\lesssim\EE^{\pi^N}\bigl[ | \EE_x[ f(Z^N(x,\xi^N) )
-f(Z_{\ell}) ] |^p \bigr]
+ \EE^{\pi^N}[|\err^N_{\star}(x,\xi^N)|^p ]\\
&&\qquad= \EE^{\pi^N}\biggl[ \biggl| F\biggl(-\frac{\ell^3}{4},\var[ Z^N(x,\xi^N)
]^{1/2}\biggr) - F\biggl(-\frac{\ell^3}{4},\var[ Z_{\ell} ]^{{1}/{2}}\biggr) \biggr|^p \biggr]
\\
&&\qquad\quad{}+ \EE^{\pi^N}[|\err^N_{\star}(x,\xi^N)|^p ]\\
&&\qquad\lesssim\EE^{\pi^N}\bigl[ |\var[ Z^N(x,\xi^N) ]^{1/2}
- \var[ Z_{\ell} ]^{{1}/{2}} |^p \bigr]
+ \EE^{\pi^N}[|\err^N_{\star}(x,\xi^N)|^p ]\\
&&\qquad\lesssim\EE^{\pi_0} \biggl| \biggl\{ \frac{\|x\|^2_{\C^N}}{N} \biggr\}
^{{1}/{2}} - 1\biggr|^p
+ \EE^{\pi^N}[|\err^N_{\star}(x,\xi^N)|^p ] \to0.
\end{eqnarray*}
In the last step we have used the fact that if $x \dist\pi_0$,
then $\frac{\|x\|^2_{\C^N}}{N} \dist\frac{\rho_1^2 + \cdots+ \rho
_N^2}{N}$ where $\rho_1, \ldots, \rho_N$
are i.i.d. Gaussian random variables $\Normal(0,1)$ so that\break
$\EE^{\pi_0} | \{ \frac{\|x\|^2_{\C^N}}{N} \}^{{1}/{2}}
- 1|^p \to0$.
\end{pf}

\begin{corollary} \label{remconcentration}
Let $p \geq1$ be a positive.
The local mean acceptance probability $\alpha^N(x)$, defined in
equation \eqref{elocaccept}, satisfies
\[
\lim_{N \to\infty} \EE^{\pi^N}[ |\alpha^N(x) - \alpha(\ell)|^p
] = 0.
\]
\end{corollary}

\begin{pf}
The function $f(z)=1 \wedge e^z$ is $1$-Lipschitz and $\alpha(\ell) =
\EE [f(Z_{\ell})]$. Also,
\[
\alpha^N(x) = \EE_x[ f(Q^N(x,\xi^N)) ] = \EE_x \bigl[f\bigl(Z^N(x,\xi
^N) + \err^N_{\star}(x,\xi^N)\bigr) \bigr]
\]
with $\err^N_{\star}(x,\xi^N) = i^N(x,\xi^N) +\err^N(x,\xi^N)$. Lemma
\ref{lemgaussianapprox} shows that\break
$\lim_{N \to\infty} \EE^{\pi^N}[\err^N_{\star}(x,\xi)^p] = 0$,
and therefore
Lemma~\ref{lemconcentration} gives the conclusion.
\end{pf}

\subsection{Drift approximation}
\label{ssecda}

This section proves that the approximate drift function $d^N\dvtx\h^s
\to\h
^s$ defined in equation \eqref{eqdrift}
converges to the drift function $\mu\dvtx\h^s \to\h^s$ of the limiting
diffusion \eqref{eqnspdemain}.

\begin{lemma}[(Drift approximation)] \label{lemdriftapprox}
Let Assumption~\ref{ass1} hold.
The drift function $d^N\dvtx\h^s \to\h^s$ converges to $\mu$ in
the sense that
\[
\lim_{N \to\infty} \EE^{\pi^N}[ \|d^N(x)-\mu(x)\|^2_s ] = 0.
\]
\end{lemma}

\begin{pf}
The approximate drift $d^N$ is given by equation \eqref{eqdrift}.
The definition of the local mean acceptance probability $\alpha^N(x)$,
given by equation \eqref{elocaccept}, show that $d^N$
can also be expressed as
\[
d^N(x) = ( \alpha^N(x) \alpha(\ell)^{-1} ) \mu^N(x)
+ \sqrt{2 \ell} h(\ell)^{-1} (\Delta t)^{-1/2} \eps^N(x),
\]
where $\mu^N(x) = -( P^N x + \C^N \nabla\Psi^N(x))$, and the
term $\eps^N(x)$ is defined by
\[
\eps^N(x)
= \EE_x[ \gamma^N(x,\xi^N) \C^{{1}/{2}} \xi^N ]
= \EE_x \bigl[ \bigl( 1 \wedge e^{Q^N(x,\xi^N)} \bigr) \C^{{1}/{2}} \xi^N \bigr].
\]
To prove Lemma~\ref{lemdriftapprox}, it suffices to verify that
\begin{eqnarray}
\lim_{N \to\infty}
\EE^{\pi^N} [ \| ( \alpha^N(x) \alpha(\ell)^{-1} )
\mu^N(x) - \mu(x) \|_s^2 ]&=& 0, \label{edriftcond1}\\
\lim_{N \to\infty} (\Delta t)^{-1}
\EE^{\pi^N} [ \|\eps^N(x) \|_s^2 ]& =& 0. \label{edriftcond2}
\end{eqnarray}
\begin{itemize}
\item
Let us first prove equation \eqref{edriftcond1}.
The triangle inequality and the Cauchy--Schwarz inequality show that
\begin{eqnarray*}
&&\bigl(\EE^{\pi^N} [ \| ( \alpha^N(x) \alpha(\ell)^{-1}
) \mu^N(x) - \mu(x) \|_s^2 ]\bigr)^2\\
&&\qquad\lesssim\EE[ |\alpha^N(x) - \alpha(\ell)|^4 ] \cdot\EE^{\pi^N}
[ \|
\mu^N(x)\|_s^4 ] \\
&&\qquad\quad{} + \EE^{\pi^N} [ \|\mu^N(x) - \mu(x)\|_s^4 ].
\end{eqnarray*}
By Remark~\ref{remts}, $\mu^N\dvtx\h^s \to\h^s$ is Lipschitz,
with a Lipschitz constant that can be chosen independent of $N$.
It follows that $\sup_{N} \EE^{\pi^N} [ \|\mu^N(x)\|_s^4 ] < \infty
$. Lemma~\ref{lemconcentration}
and Corollary~\ref{remconcentration} show that $\EE[ |\alpha^N(x) -
\alpha(\ell)|^4 ] \to0$. Therefore,
\[
\lim_{N \to\infty} \EE[ |\alpha^N(x) - \alpha(\ell)|^4 ] \cdot
\EE
^{\pi^N} [ \|\mu^N(x)\|_s^4 ] = 0.
\]
The functions $\mu^N$ and $\mu$ are globally Lipschitz on $\h^s$, with
a Lipschitz constant that can be chosen independently of
$N$, so that $\|\mu^N(x) - \mu(x)\|_s^4 \lesssim(1+\|x\|_s^4)$.
Lemma~\ref{lemmuNmu} proves that
the sequence of functions $\{\mu^N\}$ converges $\pi_0$-almost surely
to $\mu(x)$ in $\h^s$, and
Lemma~\ref{lemmchmeas} shows that\break
$\EE^{\pi^N} [ \|\mu^N(x) - \mu(x)\|_s^4 ] \lesssim\EE^{\pi_0} [
\|\mu
^N(x) - \mu(x)\|_s^4 ]$.
It thus follows from the dominated convergence theorem that
\[
\lim_{N \to\infty} \EE^{\pi^N} [ \|\mu^N(x) - \mu(x)\|_s^4 ] = 0.
\]
This concludes the proof of equation \eqref{edriftcond1}.

\item
Let us prove equation \eqref{edriftcond2}.
If the Bernoulli random variable $\gamma^N(x,\xi^N)$ were independent
from the noise term $(\C^N)^{{1}/{2}} \xi^N$,
it would follow that $\eps^N(x)=0$. In general $\gamma^N(x,\xi^N)$ is
not independent from $(\C^N)^{1/2} \xi^N$
so that $\eps^N(x)$ is not equal to zero. Nevertheless, as
quantified by Lemma~\ref{lemconcentration},
the Bernoulli random variable $\gamma^N(x,\xi^N)$ is asymptotically
independent from the current position $x$ and from the noise term $(\C
^N)^{1/2} \xi^N$.
Consequently, we can prove in equation~\eqref{eepssmall} that the
quantity $\eps^N(x)$ is small.
To this end, we establish that each component $\langle\eps(x), \pphi
_j \rangle^2_s$ satisfies
\begin{equation} \label{ecomponenteps}
\EE^{\pi^N} [ \langle\eps^N(x), \pphi_j \rangle^2_s ] \lesssim
N^{-1} \EE^{\pi^N}[\langle x,\pphi_j \rangle^2_s] + N^{-{2}/{3}} (j^s
\lambda_j)^2.
\end{equation}
Summation of equation \eqref{ecomponenteps}
over $j=1, \ldots,N$ leads to
\begin{equation} \label{eepssmall}
\EE^{\pi^N} [ \|\eps^N(x)\|_s^2 ]
\lesssim N^{-1} \EE^{\pi^N}[\|x\|_s^2 ] + N^{-
{2}/{3}} \tr_{\h^s}(\C_s)
\lesssim N^{-{2}/{3}},
\end{equation}
which gives the proof of equation \eqref{edriftcond2}.
To prove equation \eqref{ecomponenteps} for a fixed index $j \in
\mathbb{N}$, the quantity
$Q^N(x,\xi)$ is decomposed as a sum of a term, independently from $\xi
_j$, and another remaining term of small magnitude.
To this end we introduce
\begin{equation} \label{eQQi}
\cases{
Q^N(x,\xi^N) = Q^N_j(x,\xi^N) + Q^N_{j,\perp}(x,\xi^N), \vspace*{2pt}\cr
Q^N_{j}(x,\xi^N) = \displaystyle-\frac{1}{\sqrt{2}} \ell^{3/2} N^{-1/2}
\lambda_j^{-1} x_j \xi_j
- \frac{1}{2} \ell^2 N^{-2/3} \lambda_j^2 \xi_j^2\vspace*{2pt}\cr
\hspace*{62pt}{}+\err^N(x,\xi^N).}
\end{equation}
The definitions of $Z^N(x,\xi^N)$ and $i^N(x,\xi^N)$ in equations
\eqref
{eZn} and \eqref{eiN}
readily show that $Q^N_{j,\perp}(x,\xi^N)$ is independent from $\xi_j$.
The noise term satisfies $\C^{1/2} \xi^N = \sum_{j=1}^N (j^s
\lambda
_j) \xi_j \pphi_j$.
Since $Q^N_{j,\perp}(x,\xi^N)$, and $\xi_j$ are independent and $z
\mapsto1 \wedge e^z$ is $1$-Lipschitz, it follows that
\begin{eqnarray*}
\langle\eps^N(x), \pphi_j \rangle^2_s
&=& (j^s \lambda_j)^2 \bigl( \EE_x\bigl[ \bigl(1 \wedge e^{Q^N(x,\xi^N)}
\bigr) \xi_j \bigr] \bigr)^2 \\
&=& (j^s \lambda_j)^2 \bigl( \EE_x\bigl[ \bigl[\bigl(1 \wedge e^{Q^N(x,\xi
^N)} \bigr)-\bigl(1 \wedge e^{Q^N_{j,\perp}(x,\xi^N)} \bigr)\bigr] \xi_j
\bigr] \bigr)^2 \\
&\lesssim&(j^s \lambda_j)^2 \EE_x [ |Q^N(x,\xi^N)-Q^N_{j,\perp
}(x,\xi^N)|^2]\\
&= &(j^s \lambda_j)^2 \EE_x [ Q^N_{j}(x,\xi^N)^2 ].
\end{eqnarray*}
By Lemma~\ref{lemgaussianapprox} $\EE^{\pi^N}[ \err^N(x,\xi^N)^2
] \lesssim N^{-{2}/{3}}$. Therefore,
\begin{eqnarray*}
&&(j^s \lambda_j)^2 \EE^{\pi^N} [ Q^N_{j}(x,\xi^N)^2]\\
&&\qquad\lesssim(j^s \lambda_j)^2 \{ N^{-1} \lambda_j^{-2} \EE^{\pi
^N}[ x^2_j \xi^2_j ]
+ N^{-4/3} \EE^{\pi^N}[ \lambda_j^4 \xi_j^4 ]
+\EE^{\pi^N}[\err^N(x,\xi)^2 ] \} \\
&&\qquad\lesssim N^{-1} \EE^{\pi^N}[ (j^{s}x_j)^2 \xi^2_j ] + (j^s
\lambda_j)^2 (N^{-{4}/{3}} + N^{-{2}/{3}}) \\
&&\qquad\lesssim N^{-1} \EE^{\pi^N}[ \langle x,\pphi_j \rangle_s^2 ] + (j^s
\lambda_j)^2 N^{-{2}/{3}} \\
&&\qquad\lesssim N^{-1} \EE^{\pi^N}[ \langle x,\pphi_j \rangle_s^2 ] + (j^s
\lambda_j)^2 N^{-{2}/{3}},
\end{eqnarray*}
which finishes the proof of equation \eqref{ecomponenteps}.\quad\qed
\end{itemize}
\noqed\end{pf}

\subsection{Noise approximation}
\label{ssecna}

Recall definition \eqref{eqmart} of the martingale
difference $\Gamma^{k,N}.$ In this section we estimate the
error in the approximation $\Gamma^{k,N} \approx\Normal(0,\C_s)$.
To this end we introduce the covariance operator
\[
D^N(x) = \EE_x [ \Gamma^{k,N} \otimes_{\h^s} \Gamma^{k,N} |
x^{k,N}=x ].
\]
For any $x, u,v \in\h^s$, the operator $D^N(x)$ satisfies
\[
\EE[\langle\Gamma^{k,N},u \rangle_s \langle\Gamma^{k,N},v \rangle
_s |x^{k,N}=x
] = \langle u, D^N(x) v \rangle_s.
\]
The next lemma gives a quantitative version of the approximation
\mbox{$D^N(x) \approx\C_s$}.

\begin{lemma} \label{lemdiffus} Let Assumption~\ref{ass1} hold.
For any pair of indices $i,j \geq0$, the operator $D^N(x)\dvtx\h^s
\to\h
^s$ satisfies
%
%
\begin{equation}
\lim_{N \to\infty} \EE^{\pi^N} | \langle\pphi_i, D^N(x) \pphi
_j \rangle_s - \langle\pphi_i, \C_s \pphi_j \rangle_s |
= 0 \label{enoiseapprox1}
\end{equation}
and, furthermore,
%
%
\begin{equation}
\lim_{N \to\infty} \EE^{\pi^N} | \tr_{\h^s}(D^N(x)) - \tr_{\h
^s}(\C_s) |
= 0. \label{enoiseapprox2}
\end{equation}
\end{lemma}
\begin{pf}
The martingale difference $\Gamma^{N}(x,\xi)$ is given by
\begin{eqnarray} \label{emartdef}
\qquad\Gamma^{N}(x,\xi) &=&
\alpha(\ell)^{-1/2} \gamma^N(x,\xi) \C^{1/2} \xi
\nonumber
\\[-8pt]
\\[-8pt]
\nonumber
&&{}+
\frac{1}{\sqrt{2}}\alpha(\ell)^{-1/2} (\ell\Delta t)^{{1}/{2}}
\{ \gamma^N(x,\xi) \mu^N(x) - \alpha(\ell) d^N(x)\}.
\end{eqnarray}
We only prove equation \eqref{enoiseapprox2};
the proof of equation \eqref{enoiseapprox1}
is essentially identical, but easier.
Remark~\ref{remts} shows that
the functions $\mu, \mu^N\dvtx\h^s \to\h^s$ are
globally Lipschitz, and Lemma~\ref{lemdriftapprox}
shows that $\EE^{\pi^N}[ \|d^N(x)-\mu(x)\|_s^2 ] \to0$. Therefore
\begin{equation} \label{eerror1noise}
\EE^{\pi^N} [ \| \gamma^N(x,\xi) \mu^N(x)
- \alpha(\ell) d^N(x) \|_s^2 ] \lesssim1,
\end{equation}
which implies that the second term on the
right-hand side of equation \eqref{emartdef}
is $\OO(\sqrt{\Delta t})$.
Since $\tr_{\h^s}(D^N(x)) = \EE_x [ \|\Gamma^{N}(x,\xi)\|_s^2]$,
equation \eqref{eerror1noise} implies that
\[
\EE^{\pi^N}\bigl[ |\alpha(\ell) \tr_{\h^s}(D^N(x)) -
\EE_x[ \| \gamma^N(x,\xi) \C^{1/2} \xi\|_s^2 ] |
\bigr]
\lesssim(\Delta t)^{1/2}.
\]
Consequently, to prove equation \eqref{enoiseapprox2}, it suffices
to verify that
\begin{equation} \label{eerror2noise}
\lim_{N \to\infty} \EE^{\pi^N} \bigl[ \bigl| \EE_x [ \| \gamma
^N(x,\xi) \C^{1/2} \xi\|_s^2]
- \alpha(\ell) \tr_{\h^s}(\C_s) \bigr| \bigr] = 0.
\end{equation}
We have
$\EE_x [ \| \gamma^N(x,\xi) \C^{1/2} \xi\|_s^2]
= \sum_{j=1}^N (j^s \lambda_j)^2
\EE_x[ ( 1 \wedge e^{Q^N(x,\xi)} ) \xi_j^2 ]$.
Therefore, to prove equation \eqref{eerror2noise}, it suffices to establish
\begin{equation} \label{eerror3noise}
\lim_{N \to\infty} \sum_{j=1}^N (j^s \lambda_j)^2
\EE^{\pi^N} \bigl[ \bigl| \EE_x\bigl[ \bigl( 1 \wedge e^{Q^N(x,\xi)} \bigr)
\xi_j^2 \bigr] - \alpha(\ell) \bigr| \bigr] = 0.
\end{equation}
Since $\sum_{j=1}^{\infty} (j^s \lambda_j)^2 < \infty$
and $| 1 \wedge e^{Q^N(x,\xi)} | \le1$,
the dominated convergence theorem shows that \eqref{eerror3noise}
follows from
\begin{equation} \label{eerror4noise}
\lim_{N \to\infty} \EE^{\pi^N} \bigl[\bigl | \EE_x\bigl[ \bigl( 1
\wedge e^{Q^N(x,\xi)} \bigr) \xi_j^2 \bigr]
- \alpha(\ell) \bigr| \bigr] = 0\qquad \forall j \geq0.
\end{equation}
We now prove equation \eqref{eerror4noise}.
As in the proof of Lemma~\ref{lemdriftapprox}, we use the decomposition
$Q^N(x,\xi) = Q^N_j(x,\xi)+Q^N_{j,\perp}(x,\xi)$ where
$Q^N_{j,\perp}(x,\xi)$ is independent from~$\xi_j$. Therefore, since
$\operatorname{Lip}(f) = 1$,
\begin{eqnarray*}
&&\EE_x\bigl[ \bigl( 1 \wedge e^{Q^N(x,\xi)} \bigr) \xi_j^2 \bigr]\\
&&\qquad= \EE_x\bigl[\bigl ( 1 \wedge e^{Q^N_{j,\perp}(x,\xi)} \bigr) \xi_j^2 \bigr]
+ \EE_x\bigl[ \bigl[\bigl( 1 \wedge e^{Q^N(x,\xi)} \bigr)
- \bigl( 1 \wedge e^{Q^N_{j,\perp}(x,\xi)} \bigr)\bigr] \xi_j^2 \bigr]\\
&&\qquad= \EE_x\bigl[ 1 \wedge e^{Q^N_{j,\perp}(x,\xi)} \bigr]
+ \OO\bigl( \{\EE_x[ |Q^N(x,\xi)-Q^N_{j,\perp}(x,\xi)|^2]\}
^{1/2} \bigr)\\
&&\qquad= \EE_x\bigl[ 1 \wedge e^{Q^N_{j,\perp}(x,\xi)} \bigr]
+ \OO\bigl( \{\EE_x[ Q^N_j(x,\xi)^2]\}^{1/2} \bigr).
\end{eqnarray*}
Lemma~\ref{lemconcentration} ensures that, for $f(\cdot)=1
\wedge\exp(\cdot)$,
\[
\lim_{N \to\infty} \EE^{\pi^N} \bigl[ | \EE_x[ f(Q^N_{j,\perp
}(x,\xi)) ]
- \alpha(\ell) | \bigr] = 0,
\]
and the definition of $Q^N_i(x,\xi)$ readily shows that
$\lim_{N \to\infty} \EE^{\pi^N}[ Q^N_{j}(x,\xi)^2] = 0$.
This concludes the proof of equation \eqref{eerror4noise} and
thus ends the proof of Lemma~\ref{lemdiffus}.
\end{pf}
\begin{corollary} \label{corgeneralh}
More generally, for any fixed vector $h \in\h^s$, the following limit holds:
\begin{equation}
\lim_{N \to\infty} \EE^{\pi^N} | \langle h, D^N(x) h \rangle_s -
\langle h, \C_s h \rangle_s |
= 0. \label{enoiseapprox3}
\end{equation}
\end{corollary}
\begin{pf}
If $h = \pphi_i$, this is precisely the content of Proposition \ref
{propdiffapprox}.
More generally, by linearity, Proposition~\ref{propdiffapprox} shows
that this is true for
$h = \sum_{i \leq N} \alpha_i \pphi_i$, where $N \in\mathbb{N}$ is a
fixed integer.
For a general vector $h \in\h^s$, we can use the decomposition $h =
h^{*} + e^{*}$
where $h^{*} = \sum_{j \leq N} \langle h, \pphi_j \rangle_s \pphi
_j$ and
$e^{*} = h-h^*$.
It follows that
\begin{eqnarray*}
&&\bigl|\bigl( \langle h, D^N(x) h \rangle_s - \langle h, \C_s h \rangle_s \bigr)
-\bigl ( \langle h^*, D^N(x) h^* \rangle_s - \langle h^*, \C_s h^*
\rangle_s \bigr) \bigr|
\\
&&\qquad\leq|\langle h+h^*, D^N(x) (h-h^*) \rangle_s - \langle h+h^*, \C_s
(h - h^*) \rangle_s | \\
&&\qquad\leq2 \|h\|_s \cdot\|h - h^*\|_s \cdot\bigl(\tr_{\h^s}(D^N(x)) + \tr
_{\h_s}(\C_s) \bigr),
\end{eqnarray*}
where we have used the fact that for an nonnegative
self-adjoint operator $D\dvtx\h^s \to\h^s$ we have
$\langle u,Dv \rangle_s \leq\|u\|_s \cdot\|v\|_s \cdot\tr_{\h^s}(D)$.
Proposition~\ref{propdiffapprox} shows that
$\EE^{\pi^N} [ \tr_{\h^s}(D^N(x)) ] < \infty$, and
Assumption~\ref{ass1} ensures
that $\tr_{\h^s}(\C_s) < \infty$. Consequently,
\begin{eqnarray*}
&&\lim_{N \to\infty} \EE^{\pi^N} | \langle h, D^N(x) h \rangle -
\langle h, \C _s h \rangle |\\
&&\qquad\lesssim\lim_{N \to\infty} \EE^{\pi^N} | \langle h^*, D^N(x) h^*
\rangle
- \langle h^*, \C_s h^* \rangle | + \|h-h^*\|_s\\
&&\qquad= \|h-h^*\|_s.
\end{eqnarray*}
Since $\|h-h^*\|_s$ can be chosen arbitrarily small, the conclusion follows.
\end{pf}

\subsection{Martingale invariance principle} \label{secbweakconv}
This section proves that the process~$W^N$ defined in
equation \eqref{eWN} converges to a Brownian motion.
%
\begin{prop} \label{lembweakconv}%
Let Assumption~\ref{ass1} hold.
Let $z^0 \sim\pi$ and $W^N(t)$, the process defined in equation
\eqref
{eWN}, and
$x^{0,N} \dist\pi^N$, the starting position of the Markov chain
$x^N$. Then
%
%
\begin{equation}
(x^{0,N}, W^N) \longweak(z^0,W),
\end{equation}
where $\longweak$ denotes weak convergence in $\h^s \times C([0,T];\h
^s)$, and $W$ is
a $\h^s$-valued Brownian motion with covariance operator $\C_s$.
Furthermore the limiting Brownian motion $W$ is independent of the
initial condition $z^0$.
\end{prop}

\begin{pf}
As a first step, we show that $W^N$ converges weakly to $W$.
As described in~\cite{MatPillStu09}, a consequence of
Proposition $5.1$ of~\cite{Berg86} shows that in order to prove that
$W^N$ converges weakly to $W$ in $C([0,T]; \h^s)$, it suffices to
prove that
for any $t \in[0,T]$ and any pair of indices $i,j \geq0$
the following three limits hold in probability,
the third for any $\varepsilon>0$:
%
\begin{eqnarray}
\lim_{N \rightarrow\infty} \Delta t \sum_{k=1}^{k_N(T)}
\EE[\|\Gamma^{k,N}\|_s^2 | \mathcal{F}^{k,N} ]
&=& T \tr_{\h^s}(\C_s), \label{cond1}\\
\lim_{N \rightarrow\infty}\Delta t \sum_{k=1}^{k_N(t)}
\EE[ \langle\Gamma^{k,N}, \pphi_i \rangle_s \langle\Gamma^{k,N},
\pphi_j \rangle_s
| \mathcal{F}^{k,N}]
&=& t \langle\pphi_i, \C_s \pphi_j \rangle_s, \label{cond2} \\
\lim_{N \rightarrow\infty} \Delta t \sum_{k=1}^{k_N(T)}
\EE\bigl[\| \Gamma^{k,N}\|_s^2 1_{ \{\| \Gamma^{k,N} \|_s^2
\geq\Delta t \varepsilon\}} |\mathcal{F}^{k,N} \bigr] &=& 0,
\label{cond3}
\end{eqnarray}
where $k_N(t) = \lfloor\frac{t}{\Delta t} \rfloor,$ $\{\pphi_j\}$
is an
orthonormal basis of
$\h^s$ and $\FF^{k,N}$ is the natural filtration of the Markov chain
$\{
x^{k,N}\}$.
The proof follows from the estimate on
$D^N(x) = \EE[ \Gamma^{0,N} \otimes\Gamma^{0,N} | x^{0,N}=x]$
presented in Lemma~\ref{lemdiffus}.
For the sake of simplicity, we will write $\EE_k[ \cdot]$ instead
of $\EE[ \cdot|\FF^{k,N}]$. We now prove that the three conditions
are satisfied.
\begin{itemize}
\item\textit{Condition} \eqref{cond1}.
It is enough to prove that
\[
\lim\EE\Biggl| \Biggl\{ \frac{1}{\lfloor N^{{1}/{3}} \rfloor}
\sum_{k=1}^{\lfloor N^{{1}/{3}} \rfloor}
 \EE_k[ \|\Gamma^{k,N}\|_s^2 ]\Biggr\} - \tr_{\h^s}(\C_s) \Biggr| =
0,
\]
where
\[
\EE_k [ \|\Gamma^{k,N}\|_s^2 ] =
\EE_k \sum_{j=1}^N [ \langle\pphi_j, D^N(x^{k,N}) \pphi_j \rangle
_s ]
= \EE_k\tr_{\h^s}(D^N(x^{k,N})).
\]
Because the Metropolis--Hastings algorithm preserves
stationarity and $x^{0,N} \dist\pi^N$, it follows that
$x^{k,N} \dist\pi^N$ for any $k \geq0$. Therefore, for all $k \geq
0$, we have
$\tr_{\h^s}(D^N(x^{k,N})) \dist\tr_{\h^s}(D^N(x))$ where $x \dist
\pi
^N$. Consequently, the triangle inequality
shows that
\begin{eqnarray*}
&&\EE\Biggl| \Biggl\{ \frac{1}{\lfloor N^{{1}/{3}} \rfloor}
\sum_{k=1}^{\lfloor N^{{1}/{3}} \rfloor} \EE_k \|\Gamma^{k,N}\|
^2 \Biggr\} -
\tr_{\h^s}(\C_s) \Biggr|\\
&&\qquad\leq
\EE^{\pi^N} | \tr_{\h^s}( D^N(x) ) - \tr_{\h^s}(\C_s) |
\to0,
\end{eqnarray*}
where the last limit follows from Lemma~\ref{lemdiffus}.

\item\textit{Condition} \eqref{cond2}.
It is enough to prove that
\[
\lim\EE^{\pi^N} \Biggl| \Biggl\{ \frac{1}{\lfloor N^{{1}/{3}} \rfloor}
\sum
_{k=1}^{\lfloor N^{{1}/{3}} \rfloor}
\EE_k[ \langle\Gamma^{k,N},\pphi_i \rangle_s \langle\Gamma
^{k,N},\pphi_j \rangle_s
] \Biggr\} - \langle\pphi_i, \C_s \pphi_j \rangle_s \Biggr| = 0,
\]
where $\EE_k[ \langle\Gamma^{k,N},\pphi_i \rangle_s \langle
\Gamma ^{k,N},\pphi _j \rangle_s ] = \langle\pphi_i,
D^N(x^{k,N})\pphi_j \rangle_s$. Because
$x^{k,N} \dist\pi^N$, the conclusion again follows from Lemma \ref
{lemdiffus}.

\item\textit{Condition} \eqref{cond3}.
For all $k \geq1$, we have $x^{k,N} \dist\pi^N$ so that
\begin{eqnarray*}
&&\EE^{\pi^N} \Biggl| \frac{1}{\lfloor N^{{1}/{3}} \rfloor}
\sum_{k=1}^{\lfloor N^{{1}/{3}} \rfloor} \EE_{k}\bigl[\| \Gamma
^{k,N}\|_s^2
1_{ \| \Gamma^{k,N} \|_s^2 \geq N^{{1}/{3}} \eps}\bigr] \Biggr|\\
&&\qquad\leq
\EE^{\pi^N}\|\Gamma^{0,N}\|_s^2 1_{ \{ \| \Gamma^{0,N} \|_s^2 \geq
N^{{1}/{3}} \eps\}}.
\end{eqnarray*}
Equation \eqref{emartdef} shows that for any power $p \geq0$, we
have $\sup_N \EE^{\pi^N}[ \| \Gamma^{0,N} \|_s^p ] < \infty$.
Therefore the sequence $\{\|\Gamma^{0,N}\|_s^2\}$ is uniformly
integrable, which shows that
\[
\lim_{N \to\infty} \EE^{\pi^N}\|\Gamma^{0,N}\|_s^2 1_{ \{ \|
\Gamma
^{0,N} \|_s^2 \geq N^{{1}/{3}} \eps\}} = 0.
\]
\end{itemize}
The three hypothesis are satisfied, proving that $W^N$
converges weakly in $C([0,T]; \h^s)$ to a Brownian motion $W$ in $\h^s$
with covariance $\C_s$.
Therefore, Corollary $4.4$ of~\cite{MatPillStu09} shows that
the sequence $\{ (x^{0,N}, W^N) \}_{N \geq1}$
converges weakly to $(z^0, W)$ in $\h\times C([0,T], \h^s)$.
This completes the proof of Proposition~\ref{lembweakconv}.
\end{pf}

\section{Conclusion}
\label{secconclusion}

We have studied the application of the MALA
algorithm to sample from measures defined
via density with respect to a Gaussian measure
on Hilbert space. We prove that a suitably
interpolated and scaled version of the Markov chain
has a diffusion limit in infinite dimensions. There are two main conclusions
which follow from this theory: first, this work
shows that, in stationarity, the MALA algorithm
applied to an $N$-dimensional approximation of\vadjust{\goodbreak}
the target will take $\mathcal{O}(N^{1/3})$ steps
to explore the invariant measure; second, the
MALA algorithm will be optimized at an average acceptance
probability of $0.574$. We have thus significantly
extended the work~\cite{RobeRose98} which reaches similar
conclusions in the case of i.i.d. product targets.
In contrast we have considered target measures
with significant correlation, with structure motivated
by a range of applications. As a consequence our
limit theorems are in an infinite dimensional
Hilbert space, and we have employed
an approach to the derivation of the diffusion
limit which differs significantly
from that used in~\cite{RobeRose98}.
This approach was developed in~\cite{MatPillStu09}
to study diffusion limits for the RWM algorithm.

There are many possible developments of this work.
We list several of these.

\begin{itemize}

\item In~\cite{BPRSS10a} it is shown that the Hybrid
Monte Carlo algorithm (HMC) requires, for target measures
of the form \eqref{eqtarget1}, $\mathcal{O}(N^{1/4})$ steps
to explore the invariant measure. However, there is no
diffusion limit in this case. Identifying an appropriate
limit, and extending analysis to the case of target
measures \eqref{eqtarget3}, provides a challenging
avenue for exploration.

\item In the i.i.d. product case, it is known that
if the Markov chain is started ``far'' from stationarity,
a fluid limit (ODE) is observed~\cite{christensen2005scaling}.
It would be interesting to study such limits in the
present context.

\item Combining the analysis of MCMC methods for
hierarchical target measures~\cite{Beda09}
with the analysis herein
provides a challenging set of theoretical questions,
as well as having direct applicability.

\item It should also be noted that, for measures absolutely
continuous with respect to a Gaussian, there exist
new nonstandard versions of RWM~\cite{BeskStua07},
MALA~\cite{Besketal08}
and HMC~\cite{BSS10a}
for which the acceptance probability does
not degenerate to zero as dimension $N$ increases.
These methods may be expensive to implement when
the Karhunen--Lo\`eve basis is not known explicitly,
and comparing their overall efficiency with that
of standard RWM, MALA and HMC is an interesting
area for further study.

\item It is natural to ask whether
analysis similar to that undertaken here could be developed
for Metropolis--Hastings methods applied to other
reference measures with a non-Gaussian product structure.
Particularly, the Besov priors of~\cite{las09} provide
an interesting class of such reference measures, and the
paper~\cite{DHS11} provides a machinery for analyzing
change of measure from the Besov prior, analogous to that
used here in the Gaussian case.
Another interesting class of reference measures
are those used in the study of uncertainty quantification
for elliptic PDEs: these have the form of an
infinite product of compactly supported uniform distributions;
see~\cite{SS11}.
\end{itemize}

\section*{Acknowledgments}
Part of this work was done when A.~H. Thi{\'e}ry  was visiting the Department of
Statistics at
Harvard University, and we thank this institution for its
hospitality.We also thank the referee for his/her very useful
comments.

%

%


\printaddresses

\end{document}